\documentclass[11pt]{article}

\usepackage{psfig}
\usepackage{amssymb,amsmath}

\author{Fran\c{c}ois Ledrappier$^{(1)}$, Michael Shub$^{(2)}$, Carles
Sim\'o$^{(3)}$, Amie Wilkinson$^{(4)}$}

\title{Random versus deterministic exponents in a rich family of
diffeomorphisms}

\date{}

\newtheorem{theorem}{Theorem}[section]
\newtheorem{prop}[theorem]{Proposition}
\newtheorem{lemma}[theorem]{Lemma}
\newtheorem{claim}[theorem]{Claim}
\newtheorem{cor}[theorem]{Corollary}
\newtheorem{quest}[theorem]{Question}
\newtheorem{definition}{Definition}[section]
\newtheorem{remark}[theorem]{Remark}

\def\eproof{\;{\bf QED}}

\def\eps{\varepsilon}
\def\bar{\overline}
\def\dim{\hbox{dim}}
\def\exp{\hbox{exp}}
\def\Diff{\hbox{Diff$\mskip 1.5 mu$}}
\def\GL{GL(n,{\bf C})}
\def\GLR{GL(n,{\bf R})}
\def\GLDR{GL(2,{\bf R})}
\def\GLDPR{GL^+(2,{\bf R})}
\def\GLDMR{GL^-(2,{\bf R})}
\def\GR{G_{n,k}({\bf C})}
\def\Cn{{\bf C}^n}
\def\C{{\bf C}}
\def\Un{U(n,{\bf C})}
\def\SO2{SO(2,{\bf R})}
\def\O2{O(2,{\bf R})}
\def\RE{{\mathbb R}}
\def\S1{{S}^1}
\def\CF{{\cal F}}
\def\CM{{\cal M}}
\def\F{\mathcal F}
\def\H{\mathcal H}

\newcommand{\Nset}{\bf N }
\newcommand{\Qset}{\bf Q}
\newcommand{\Rset}{\bf R}
\newcommand{\Sset}{S}
\newcommand{\Cset}{\bf C}
\newcommand{\dd}{\,d}

\begin{document}

\maketitle
 
\bigskip
\bigskip        

\centerline{\em Dedicated to Yakov Sinai on his 65th birthday.}

\begin{abstract} We study, both numerically and theoretically,
the relationship between the
random Lyapunov exponent of a family of area
preserving diffeomorphisms of the $2$-sphere and the mean of the
Lyapunov exponents of the individual members.
The motivation for this study is the hope that a rich enough family of
diffeomorphisms will always have members with positive Lyapunov
exponents, that is to say, positive entropy. At question is what
sort of notion of richness would make such a conclusion valid.
One type of richness of a family -- invariance under the left action
of $SO(n+1)$ -- occurs naturally in the context of volume preserving
diffeomorphisms of the $n$-sphere. Based on some
positive results for families linear maps obtained by
Dedieu and Shub, we investigate the exponents of such a family
on the $2$-sphere.  Again motivated by the linear case, we
investigate whether there is in fact a lower
bound for the mean of the Lyapunov exponents in terms of the
random exponents (with respect to  the push-forward of Haar measure
on $SO(3)$) in such a family.  The family ${\mathcal F}_\eps$ that
we study contains
a twist map with stretching parameter $\eps$.

In the family ${\mathcal F}_\eps$, we find
strong numerical evidence for
the existence of such a lower bound on mean Lyapunov exponents,
when the values of the
stretching parameter $\eps$ are not too small. Even moderate values of
$\eps$ like $\eps\ge 10$ are enough to have an
average of the metric entropy larger than that of the
random map. For small $\eps$ the estimated average entropy seems positive but
is definitely much less than the one of the random map. The numerical
evidence is in favor of the existence of exponentially small lower and
upper bounds (in the present example, with an analytic family).

Finally, the effect of a small randomization of fixed size $\delta$ of the
individual elements of the family ${\mathcal F}_\eps$ is considered. Now
the mean of the local random exponents of the family is indeed 
asymptotic to the random exponent of the entire family as $\eps$ 
tends to infinity.

\end{abstract}

\vfill

{\small{
\noindent $^{(1)}$ Centre de Math\'ematiques, \'Ecole Polytechnique, 91128
Palaiseau Cedex, France\\
\noindent $^{(2)}$ IBM T. J. Watson Research Center,  P. O. Box 218,
Yorktown Heights, NY 10598, USA\\
\noindent $^{(3)}$ Dept. de Matem\`atica. Aplicada i An\`alisi, Univ. de
Barcelona, 08071 Barcelona, Spain\\
\noindent $^{(4)}$ Mathematics Department, Northwestern University, Evanston IL
60208-273, USA\\
\noindent E-mail: \texttt{ledrappi@math.polytechnique.fr, mshub@us.ibm.com,
carles@maia.ub.es,}\\
\noindent \phantom{E-mail:}$\,$ \texttt{wilkinso@math.northwestern.edu}
}}

\eject

\section{Introduction} \label{S=intro}

Numerical experiments with area-preserving surface diffeomorphisms often produce
the following dynamical picture: elliptical islands floating in ergodic seas.
A reasonable guess is that these ergodic seas typically have positive measure,
and further, that the Lyapunov exponents on these seas are on average nonzero.
An example of tiny elliptical islands in the context of differential
equations can be found in \cite{NSS}, where all rough numeric tests are
in favor of ergodic behavior.

In this paper, we add to the pile of experimental evidence
in favor of this conjecture.\footnote{\bf Due to space limitations,
this version of
the paper does not contain figures.  A version with the
figures is available at either
http$:\slash\slash$www.math.northwestern.edu$\slash_{\,\widetilde{}}\,$wilkinso$\slash$papers ($8\frac{1}{2}\times 11$ format)
or 
http$:\slash\slash$www.maia.ub.es$\slash$dsg$\slash$2002$\slash$index.html (A4 format). }
We also discuss a {\em possible}
theoretical approach to finding positive Lyapunov exponents in certain
families of area-preserving diffeomorphisms of the sphere $\Sset^2$.
The possibility of such an approach was discussed in \cite{BPSW}.
The families we consider are not obtained from a specific set of
equations, but from the following construction.  Let $f: \Sset^2\to \Sset^2$
be an area-preserving diffeomorphism of the round sphere,
and let $SO(3)$ be the isometry group of $\Sset^2$.
Let
$${\mathcal F} = \{g\circ f\,\mid\, g\in SO(3)\}$$
be the left $SO(3)$-coset of $f$ in $\Diff(\Sset^2)$, and
let $\nu$ be the push-forward of Haar measure on $SO(3)$
to ${\mathcal F}$.
Provided that $f$ is not itself an isometry, the family ${\mathcal F}$
has nonzero {\em random} Lyapunov exponents with respect to $\nu$
(see Proposition~\ref{P=random} below).
The question this paper addresses is whether these random exponents
can somehow be connected to the Lyapunov exponents
of individual members of ${\mathcal F}$, at least {\em on $\nu$-average}.

To test whether there might be such a connection, we chose
$f$ to be a twist map, all of whose Lyapunov exponents are zero.
The resulting family ${\mathcal F}$ has similarities
to the standard family on the $2$-torus. The dynamics of
the individual elements of ${\mathcal F}$ and how they depend on parameters
is an interesting topic, but we only study here some key properties in the case
of small $\eps$.  We mainly focus on two quantities, the
{\em random exponent} $R(\nu)$ and the {\em average exponent}
$\Lambda(\nu)$, which we now define.

Let $\mu$ be Lebesgue measure on $\Sset^2$  normalized to be a
probability measure.
Suppose for now that $\nu$ is an arbitrary Borel probability
measure supported on a subset ${\mathcal F}$ of $\Diff_{\mu}(\Sset^2)$, the
space of $\mu$-preserving diffeomorphisms of $\Sset^2$.
For $f\in {\mathcal F}$, the largest Lyapunov of $f$
at $x\in \Sset^2$ is found by  computing the limit:
\begin{equation}
\lim_{n\to\infty} {1 \over n} \log \|T_xf^n\| = \lambda_1(x,f),
\label{e=lya1}
\end{equation}
which exists for $\mu$-almost every $x$ by the
subadditive ergodic theorem.
We define the {\em average exponent of $f$}
to be
\begin{equation}
\lambda(f) = \int_{\Sset^2}\lambda_1(f,x)\dd\mu(x), \label{e=lya2}
\end{equation}
and the {\em average exponent of $\nu$} to be
\begin{equation}
\Lambda(\nu) = \int_{\Diff^r_{\mu}(\Sset^2)} \lambda(f)\dd\nu(f).
\label{e=lya3}
\end{equation}

Rather than iterate a single diffeomorphism $f\in \CF$, we might choose
instead a sequence of diffeomorphisms $\{f_1, f_2, \ldots\}\subset
{\mathcal F}$ and form their composition:
$$f^{(n)} := f_n\circ f_{n-1}\circ \cdots \circ f_1.$$
If the sequence is chosen to be independent and identically distributed
with respect to $\nu$, then almost surely the limit
\begin{equation}
\lim_{n\to\infty} \frac{1}{n}\log \|T_xf^{(n)}\|=:R(x,(f_i)_1^\infty,\nu)
\label{e=ran1}\end{equation}
will exist, for $\mu$-almost every $x$.
(This too follows from the subadditive ergodic theorem, applied in the
appropriate context). Further, the integral of $R(x,(f_i)_1^\infty,\nu)$
with respect to $\mu$ is almost surely independent of the sequence
$(f_i)_1^\infty$. We define the {\em random exponent} of ${\nu}$ to be
this integral:
\begin{equation}
R(\nu) = \int_{\Sset^2} R(x,(f_i)_1^\infty,  \nu) \dd\mu(x).
\label{e=ran2}\end{equation}
(see Kifer for an introduction to the subject of random diffeomorphisms
and their exponents. We also give a self-contained introduction
in section~\ref{S=Facts}).  The random exponent $R(\nu)$ is usually positive,
unless $\nu$ is fairly degenerate \cite{Car}.

The quantity $\Lambda(\nu)$ is mysterious from a computational
perspective, but useful from a dynamical one. The quantity
$R(\nu)$ is relatively easy to estimate and is often positive.

Our goal is to understand in general if there is a notion of
richness for a probability measure $\nu$ on the volume preserving
diffeomorphisms of a closed manifold $M$ such that the positivity
of $R(\nu)$ implies the positivity of $\Lambda(\nu)$. Here we are
investigating whether the $SO(3)$ invariance of the measure $\nu$
on $\Diff_{\mu}(\Sset^2)$ might provide such a notion of richness.
In \cite {BPSW} we asked if even more might hold, that we might
bound $\Lambda(\nu)$ from below in terms of $R(\nu).$

\begin{quest} \label{Q6}
 Is there a positive constant $C$ --- perhaps $1$ ---
 such that  $\Lambda(\nu) \geq C R(\nu)$?
\end{quest}

Some motivation for Question \ref{Q6} can be found in similar question
for the iterates of linear maps (see \cite{DeSh} where an affirmative
answer to the analogue of Question \ref{Q6} is proven with $C=1$).
In Section~\ref{S=Heuristics}, we describe a theoretical framework in which
to address Question~\ref{Q6} and related questions.
We discuss the linear case in Section~\ref{S=linear}.

Returning to the specific family of diffeomorphisms mentioned earlier,
we now describe the experiment in more detail.

For $\eps>0$, we define a one-parameter family of twist maps $f_\eps$ as
follows. Express $\Sset^2$ as the sphere of radius $1/2$ centered at $(0,0)$
in $\Rset\times \Cset$, so that the coordinates $(r,z)\in \Sset^2$ satisfy
the equation
$$|r|^2 +|z|^2 = 1/4.$$
In these coordinates define a twist map $f_\varepsilon: \Sset^2\to \Sset^2$,
for $\varepsilon>0$, by
$$f_\varepsilon(r,z) = (r, \exp(2\pi i (r+ 1/2) \varepsilon) z).$$
Let $\F_\varepsilon$ be the orbit  $SO(3)f_\varepsilon$.  Let $\nu$ be the
push-forward of Haar measure on $SO(3)$. We denote the resulting random and
average Lyapunov exponents by $R(\eps)$ and $\Lambda(\eps)$, respectively.

The numerical results are described in Section~\ref{S=Numerics}.
To summarize these results, it appears that the inequality
$\Lambda(\eps)\geq R(\eps)$ is satisfied for large $\eps$, and it
is definitely not satisfied for small $\eps$. We now have rigorous
results to confirm some of these observations. The strongest of
these results is proved in Section~\ref{S=smalleps}: for $\eps$
close to $0$, there is no $C>0$ satisfying the inequality in
Question~\ref{Q6}: in fact, we show in Corollary~\ref{C=corol}
that for small $\eps$, $\Lambda(\eps)$ is less than $\eps^3 \sim
R(\eps)^{3/2}$.  The numerics support an upper  bound on
$\Lambda(\eps)$ that is exponentially small, and we show in
Theorem~\ref{T=final} that on most of ${\mathcal F}_\eps$ this is
indeed the case.
On the other hand, $\Lambda(\eps)$ does appear to be positive for
positive $\eps$ in all of the range where we can meaningfully compute.
Section~\ref{S=smalleps} also contains a study of the bifurcation
structure of fixed points inside the family ${\mathcal F}_\eps$.
We include this analysis because it sheds light on where, in both
$\F_\eps$ and $\Sset^2$, new elliptic periodic points and their
surrounding islands are produced. Homoclinic bifurcations also
give rise to horseshoes, which are associated (at least
heuristically) with positive measure sets with nonzero exponents.

For the case of large $\eps$, we show in Section~\ref{S=epstoinfty},
that an inequality like that in Question~\ref{Q6} is satisfied
when a small amount of noise is introduced.
More precisely, we first prove in Section~\ref{S=Heuristics} some results
about the quantities  $R(\eps)$ and $\Lambda(\eps)$ and a third quantity
$R(\eps,\delta)$, which measures the exponents of the ``in-between''
process in which each element of $\F_\eps$ has added noise in a $\delta$-ball
inside $\F_\eps$ (see Section \ref{S=Heuristics} for details). In particular,
we prove that any element $h$ of a family $\F=SO(3)f$ described above will
have average  $\delta$-diffused exponents $R(h,\nu,\delta)$ that are positive,
unless $f$ is an isometry. In addition, there exists a stationary measure
$m_{h,\delta}$ for such a process on the projective bundle
$P\Sset^2$ that is absolutely continuous
with respect to Liouville measure and projects to Lebesgue measure
$\mu$ on $\Sset^2$; this measure is unique among stationary measures
with these properties and is the unique fixed point of
a ``simple'' linear operator. The integrated measure
$$m_{\delta} = \int_{\F}  m_{h,\delta} d\nu(h)$$
determines $R(\nu,\delta)$, which is the average of the
$\delta$-diffused exponents $R(h,\nu,\delta)$ over $h\in\F$.
Whenever $m_{\delta}$ is equal to Lebesgue measure $m$, we have the equality:
$R(\nu,\delta)=\Lambda(\nu)$. For the family $\F_\eps$, denote by
$m_{\eps,\delta}$ this integrated measure. In Section~\ref{S=epstoinfty},
we prove (Theorem~\ref{t=converge}) that for any $\delta>0$,
$\lim_{\eps\to\infty} m_{\eps,\delta} = m$.  Using this result,
we prove that if enough noise is introduced, then the inequality in
Question~\ref{Q6} is satisfied as $\eps\to\infty$; in particular, we show in
Corollary~\ref{c=inftyineq} that $R(\eps,\delta) - R(\eps)$
tends to $0$ as $\delta\to 0$ and $\eps\to\infty$ sufficiently
quickly, for instance $\eps > \delta^{-25}$.  Results of a similar
nature for the standard family were
obtained by Carleson and Spencer \cite{CaSp} and are described
in Section~\ref{S=epstoinfty} below.

We suspect that a further study of these  measures
$m_{\delta}$ would be interesting.   Even for an
$SO(n)$- or $SU(n)$-invariant family of matrices, the
properties of the analogous ``in-between'' measures $m_\delta$ are,
for the most part, unknown.
In Section~\ref{S=linear}, we discuss what is known about
these measures.  For $SO(2)$, we prove that
$m_{\delta}$ is Lebesgue measure for all $\delta>0$,
and for $SU(n)$, $m_{\delta}$ is {\em not} Lebesgue measure
if $\delta$ is sufficiently small.

We thank Roy Adler, Victor Klepstyn, Yuri Kifer, Marco Martens,
Lai-Sang Young and especially Tolya Katok for conversations about
this work. The computing facilities of the Dynamical Systems Group
of the University of Barcelona have been widely used. We thank the
supporting institutions.
Michael Shub was partly supported by NSF Grant $\#$DMS-9988809.
Carles Sim{\'o} was partly supported
by grants DGICYT BFM2000-805, CIRIT 2001SGR-70 and INTAS00-221.
Amie Wilkinson was  partly supported by NSF Grant $\#$DMS-0100314.

\section{Background on random transformations and exponents}\label{S=Facts}

In this section we introduce some notation and gather together
some facts and propositions. What we have to say in this
section and the next is standard and can be found
for example in \cite {ki1}, \cite {ki2}, \cite {ki3}, \cite{L-Q}, \cite
{Car},
\cite {G-M} in most cases in greater
generality. We have outlined proofs here in order to be reasonably self
contained.

If $\H\subset \Diff(M)$, and $\nu$ is a probability measure on $\H$,
then $(\H,\nu)$ generates a random process given by selecting
an independent, $\nu$-distributed sequence $(h_i)_1^\infty \subset \H$ and
forming the compositions:
$$h^{(n)} = h_n\circ h_{n-1}\circ  \cdots \circ h_1.$$
To study all possible outcomes of this experiment, we introduce
the following auxiliary spaces and transformations: the
shift space, $\H^\infty :=\Pi_{j=1}^\infty \H$,
the one-sided shift $\sigma: \H^\infty\hookleftarrow$ given by:
$$\sigma(h_1,h_2,\ldots) = (h_2,h_3,\ldots),$$
and the skew product $\tau: \H^\infty \times M \hookleftarrow$
given by:
$$\tau( (h_i)_1^\infty, x) = (\sigma((h_i)_1^\infty), h_1(x)).$$
Then $\sigma$ has a natural invariant measure $\nu^\infty$,
the product measure induced by $\nu$, but {\em a priori} $\tau$ has no
preferred invariant measure.

\begin{definition} Let $\nu$ be a probability measure on $\H\subset
\Diff(M)$.
A measure $\mu$ on $M$ is {\em stationary} for the random process given by
$(\H,\nu)$ if any of the following equivalent conditions is satisfied:
\begin{enumerate}
\item $\tau_\ast(\nu^\infty\times\mu) = \nu^\infty\times\mu$
\item $ev_*(\nu\times\mu) = \mu$,
where $ev: \H\times M\to M$ is the evaluation map:
$$ev(h,x) = h(x)$$
\item $\mu\star \nu = \mu$, where $\star$ is the convolution
operator defined by:
$$\mu\star\nu (A) = \int_{\H} \mu(h^{-1}(A))\dd\nu(h),$$
for every $\mu$-measurable subset $A\subset M$.
\end{enumerate}
\end{definition}

Stationary measures always exist \cite{ki1} and are the random analogue
of invariant measures in the nonrandom setting.  Part of the focus of this
paper is to find natural stationary measures in the case where $M =
T_1\Sset^2$
and $\H$ and $\nu$ are derived from Haar measure on $SO(3)$.

Given an injective linear map $A:V
\rightarrow W$ between normed vector spaces we denote by
$A_\sharp$ the induced map from the unit sphere in $V$ to the
unit sphere in $W$, which is defined by
$v\rightarrow\frac{A(v)}{\|A(v)\|}.$ We use
the same notation for the induced map on the
projective space $PV$.  We denote the tangent bundle
of $\Sset^2$ by $T\Sset^2$, the unit tangent bundle by $T_1\Sset^2$,
and the projective bundle by $P\Sset^2$. We
let
$m$ denote the normalized Liouville measure on $P\Sset^2$,
so $m$ is a probability measure
which pushes forward under projection to $\Sset^2$ to $\mu$. The
fibers of $T\Sset^2$ and  $P\Sset^2$ over a point $z \in \Sset^2$ are
denoted by $T_z\Sset^2$ and $P_z\Sset^2$. For any manifolds $M,N$ and
differentiable map $F:M \rightarrow N$ the derivative of $F$ at $x \in M$
is denoted by
$T_xF$; for $v\in T_xM$ we will usually
write ``$TF v$'' instead of $T_x F(v)$.
Finally, we denote by $F_\sharp: T_1M\to T_1N$ the
map that covers $F$ and is $(T_xF)_\sharp$ on the fiber over $x\in M$.
Since the tangent map to $g\in SO(3)$ preserves unit tangent vectors,
we will write ``$g$'' for $g_\sharp$.

Now let $f\in\Diff_\mu(\Sset^2)$, let ${\mathcal F} = \{g\circ f\,\mid\,
g\in SO(3)\}$,
and let $\nu$ be the push-forward to ${\mathcal F}$ of Haar measure on
$SO(3)$.
Let $m$ be normalized Liouville measure on $P\Sset^2$.
Associated to ${\mathcal F}$ we then
have the set
$${\mathcal F_{\sharp}} = \{h_\sharp|h \in \mathcal F\}=\{g \circ
f_\sharp\,\mid\, g\in SO(3)\},$$
and the measure $\nu_\sharp$, the push-forward to $\Diff(P\Sset^2)$ of
Haar measure on $SO(3)$.
Let $\sigma:\F^\infty\hookleftarrow$, $\tau:\F^\infty\times \Sset^2
\hookleftarrow$,
$\sigma_\sharp: {\F_\sharp}^\infty\hookleftarrow$, and
$\tau_\sharp: {\F_\sharp}^\infty\times P\Sset^2\hookleftarrow$
be the associated auxiliary transformations
to the random processes generated by $(\F,\nu)$ and
$(\F_\sharp, \nu_\sharp)$
respectively.

\begin{lemma}
The measures $\mu$ and $m$ are stationary for $\nu$ and $\nu_\sharp$
respectively.

The transformations $\tau$, $\sigma$  and
$\sigma_\sharp$ are ergodic with respect to  $\nu^\infty\times \mu$,
$\nu^\infty$, and $\nu_\sharp^\infty$, respectively.
\end{lemma}

\begin{proof} It is straightforward to check that they are stationary.

Ergodicity is not much harder to check.
\eproof\end{proof}

\bigskip
Now we can compute $R(\nu)$ more explicitly:
\begin{prop}\label{P=random}
Let $f\in\Diff_\mu(\Sset^2)$, let ${\mathcal F} = \{g\circ f\,\mid\, g\in
SO(3)\}$,
and let $\nu$ be the push-forward to ${\mathcal F}$ of Haar measure on
$SO(3)$.
Let $m$ be normalized Liouville measure on $P\Sset^2$, the projective
bundle of $\Sset^2$.
Then
$$R(\nu) = \int_{P\Sset^2} \log\|Tf v\|\dd m(v).$$
Moreover, $R(\nu)>0$, unless $f$ is an isometry.
\end{prop}
\begin{proof}
We first apply Birkhoff's Ergodic Theorem to the measure-preserving
transformation $\tau_\sharp:\F_\sharp^\infty\times
P\Sset^2\hookleftarrow$
and the function
$$\psi((h_{i\sharp})_1^\infty, v)= \log\|Tf v\|$$
to obtain that
\begin{eqnarray*}
\lim_{p\to\infty} \frac{1}{p} \log\|Th_{p} \dots Th_{1} (v)\|
&= &\lim_{p\to \infty} \frac{1}{p} \sum_{j=1}^{p}\log \|Tf(h_{j\sharp}
\cdots h_{1\sharp}(v))\|\\
& = & \lim_{p\to \infty} \frac{1}{p}
\sum_{j=0}^{p-1}\psi(\tau_\sharp^j((h_{i\sharp})_1^\infty, v))\\
& =: & L((h_{i\sharp})_1^\infty, v)
\end{eqnarray*}
exists a.e. in  $\F_\sharp^\infty\times P\Sset^2$.
The integral of this limit $L$ is
\begin{eqnarray*}
\int L((h_{i\sharp})_1^\infty, v)
\dd (\nu_\sharp^\infty\times m) &=& \int_{\F_\sharp^\infty\times
P\Sset^2} \psi\dd(\nu_\sharp^\infty \times m)  \\
&=&\int_{P\Sset^2} \log\|Tf v\| \dd m(v).
\end{eqnarray*}

Next, we apply Oseledec's theorem to the
map $\tau:\F^\infty\times \Sset^2\hookleftarrow$ and the cocycle
$((h_{i})_1^\infty, x)\mapsto T_x h_1$.
We obtain that for almost all $(h_i)_1^\infty$, almost all $x\in \Sset^2$,
and for almost all $v\in T_{1,x}\Sset^2$,  the limit
$$K((h_i)_1^\infty, x):= \lim_{p\to\infty} \frac{1}{p} \log \|  Th_p
\ldots Th_1(v)\|$$
exists, is independent of $v$, and
has $\nu^\infty\times \mu$-integral equal to $R(\nu)$.
The function:
\begin{eqnarray*}
K((h_{i})_1^\infty) &:=&  \int_{\Sset^2} K((h_{i})_1^\infty, x)\dd\mu(x)\\
&=& \int_{x\in \Sset^2} \int_{u\in T_{1,x}\Sset^2}  \lim_{p\to\infty}
\frac{1}{p}\log \|  Th_p \ldots T_xh_1(u)\| \dd u\dd\mu(x)\\
&=&  \int_{P\Sset^2}\lim_{p\to\infty} \frac{1}{p}\log \|  Th_p \ldots
Th_1(v)\| \dd m(v)
\end{eqnarray*}
is $\sigma$-invariant and has $\nu^\infty$-integral equal to
$R(\nu)$; ergodicity of $\sigma$ implies that
it is a.e.  constant and therefore equal to $R(\nu)$.
We conclude that
\begin{eqnarray*}
R(\nu) &=& \int_{P\Sset^2}
\log \|Tf(v)\| \dd m(v)
\end{eqnarray*}

It remains to see that $\int_{P\Sset^2}
\log \|Tf v\|\,dm(v)> 0$ if $f$ is not an isometry.
For this we use the following elementary lemma.
\begin{lemma} Let $A\in SL(2,\Rset)$. Then
the Jacobian of $A_\sharp$ with respect to Lebesgue measure
on $\Sset^1$ is given by:
$$\hbox{Jac}(A_\sharp)(v) = \|A v\|^{-2},$$
for $v\in \Sset^1$.
\end{lemma}

Since $f$ preserves $\mu$, it follows from this lemma
that the Jacobian of $f_\sharp$ with respect to $m$ at
$v \in T_{1}\Sset^2$ is $\|Tf (v)\|^{-2}$.
Since $Tf_\sharp$ is a diffeomorphism,
\begin{eqnarray*}
\int_{P\Sset^2}\|Tf (v)\|^{-2}\dd m(v) &=&
\int_{P\Sset^2} \hbox{Jac}(f_\sharp)(v)\\
&=&1.
\end{eqnarray*}
By Jensen's inequality,
\begin{eqnarray*}
\int_{P\Sset^2}\log \|Tf (v)\|^{-2}\dd m(v) &\leq &
\log(\int_{P\Sset^2} \|Tf (v)\|^{-2}\dd m(v))\\
&=& 0
\end{eqnarray*}
with inequality holding unless
$\log \|Tf v\|$ is constant and equal to $0$.  Rearranging
the inequality, we see that,
unless $f$ is an isometry, we must have
$\int_{P\Sset^2}\log \|Tf (v)\|\dd m(v) >  0$.
\eproof\end{proof}

\section{A theoretical framework} \label{S=Heuristics}

\subsection{Connecting $R(\eps)$ to $\Lambda(\eps)$}

In this section, we attempt to interpolate between $R(\eps)$ and
$\Lambda(\eps)$ via a third quantity, $R(\eps, \delta)$,
which we call the {\em random $\delta$-diffused exponent}.
When $\delta$ is greater than or equal to the radius
of $SO(3)$, $R(\eps, \delta)$ is
equal to $R(\eps)$; as $\delta$ approaches $0$, $R(\eps,\delta)$
approaches (in limsup) a lower bound for $\Lambda(\eps)$.
Roughly speaking,  $R(\eps, \delta)$ is the exponent
(averaged over $\F_\eps$) obtained by introducing random perturbations
(viewed as noise) of order
$\delta$ to each element of $\F_\eps$, staying within
the family $\F_\eps$.
In Lemma~\ref{l=limit} we show that
$\limsup_{\delta\to 0} R(\eps,\delta)\leq \Lambda(\eps)$.
On the other hand, we derive in Proposition~\ref{p=formula}
a formula for $R(\eps,\delta)$:
\begin{eqnarray*}
R(\eps,\delta) & = & \int_{P\Sset^2} \log\|Tf v\| \dd
m_{\eps,\delta}(v).
\end{eqnarray*}
The probability measure $m_{\eps,\delta}$ in this formula has nice
properties:
it projects to Lebesgue measure $\mu$ on
$\Sset^2$, and is absolutely continuous with respect to $m$,
with smooth density.

Now, recall (Proposition~\ref{P=random}) that
\begin{eqnarray*}
R(\eps) & = & \int_{P\Sset^2} \log\|Tf v\| \dd m(v).
\end{eqnarray*}
If it were the case that $m_{\eps,\delta}\to m$ as $\delta\to 0$,
then it would follow that:
\begin{eqnarray*}
\Lambda(\eps) & \geq & \limsup_{\delta\to 0} R(\eps,\delta)\\
& = &  \limsup_{\delta\to 0} \int_{P\Sset^2} \log\|Tf v\| \dd
m_{\eps,\delta}(v)\\
& = & \int_{P\Sset^2} \log\|Tf v\| \dd m(v)\\
& = & R(\eps).
\end{eqnarray*}
Hence the properties of this measure $m_{\eps,\delta}$ are potentially
quite
interesting with regard to Question~\ref{Q6}.

Here we collect some properties of  $m_{\eps,\delta}$
and $R(\eps,\delta)$.
First of all, $R(\eps,\delta)$ is
always positive for $\delta > 0$ (in fact, we prove
in Corollary~\ref{c=formula} that this is true not just on
average, but for individual elements of $\F_\eps$).
In other words, introducing noise (no matter how small)
to an element $h\in \F_\eps$ invariably produces positive exponents.

The measure $m_{\eps,\delta}$ has additional properties as well.
We prove  that we can write:
$$m_{\eps,\delta} = \int_{SO(3)}  m_{g,\eps,\delta} \,dg,$$
where, for each $g\in SO(3)$, $m_{g,\eps,\delta}$ is
the unique probability measure on $P\Sset^2$ with the properties:
\begin{enumerate}
\item $m_{g,\eps, \delta}$ is stationary for the $\delta$-diffused
process about $gf_\eps$;
\item $m_{g,\eps, \delta}$ is absolutely continuous with respect to $m$,
with smooth density;
\item $m_{g,\eps,\delta}$ projects to Lebesgue measure $\mu$ on $\Sset^2$.
\end{enumerate}
As part of the proof, we
show that each measure $m_{g,\eps,\delta}$ is the unique fixed point
of a ``simple'' linear operator.

Finally, from the way $m_{\eps,\delta}$ is constructed,
it follows that  $m_{\eps,\delta}$
shares all of the symmetries of $f_\eps$.  In particular,
the density $\varphi_{\eps,\delta}$ is invariant under all rotations
that fix the North pole.  The further study of these measures
$m_{g,\eps,\delta}$
and $m_{\eps,\delta}$ might be of independent interest.
We discuss the linear version of $m_{\eps,\delta}$ in Section~\ref{S=linear}.
In Section~\ref{S=epstoinfty}, we examine the behavior of $m_{\eps,\delta}$
as $\eps\to\infty$.

We now turn to the proofs of  assertions 1.-3..  We first prove
a standard semicontinuity result for random exponents.

\begin{lemma}\label{l=expolimit}  Let $\{\gamma_i\}$ be a sequence of
probability measures on $\H\subset\Diff(M)$ that converges
weakly to a probability measure $\gamma$.
Suppose that $\mu$ is stationary for the
random process generated by $(\H,\gamma_i)$, for every $i$,
(for example, if $\H\subset\Diff_\mu(M)$).
Then
$$\limsup_{\gamma_i \rightarrow \gamma} R(\gamma_i) \leq R(\gamma).$$
\end{lemma}

\begin{proof}
Let $a_n((h_i)_1^\infty,x) = \log\|T_xh^{(n)}\|$.  Then
$a_n: \H^\infty \times M\to {\bf R}$ is subadditive with respect
to $\tau$. By the subadditive ergodic theorem  it then follows that
\begin{eqnarray*}
R(\gamma)   & = & \int_{\H^\infty\times M}\lim_{n\to\infty} \frac{1}{n}a_n
\dd(\gamma^\infty\times\mu)\\
&=&\lim_{n\to\infty} \frac{1}{n} \int_{\H^\infty\times M} a_n
\dd(\gamma^\infty\times\mu)\\
&=&\inf_n\frac{1}{n}\int_{\H^\infty\times M}
a_n\dd(\gamma^\infty\times\mu).
\end{eqnarray*}

Now for any fixed $n$ we have
\begin{eqnarray*}
\limsup_{\gamma_i \rightarrow \gamma} R(\gamma_i) &\leq&\limsup_{\gamma_i
\rightarrow \gamma} \frac{1}{n}\int_{\H^\infty\times M}
a_n\dd(\mu\times\gamma_i^\infty)\\
&=&\frac{1}{n}\int_{\H^\infty\times M} a_n\dd(\mu\times\gamma^\infty).
\end{eqnarray*}
So
\begin{eqnarray*}
\limsup_{\gamma_i \rightarrow \gamma} R(\gamma_i)&\leq& \inf_n
\frac{1}{n}\int_{\H^\infty\times M} a_n\dd(\mu\times\gamma^\infty)\\
&=&R(\gamma).
\end{eqnarray*}
$\!\!$\eproof \end{proof}

We will apply this lemma to the situation where $\gamma_i$
is supported on a small ball in $\F_\eps$ converging,
as $i\to\infty$, to a Dirac measure supported on an
element of $\F_\eps$.

Let $\delta>0$, and let $U_\delta$ be a
symmetric $\delta$ ball around the identity in $SO(3)$. Give
$U_{\delta}$ the restriction of Haar measure, normalized to be a
probability measure and similarly for
$\F_{g,\eps,\delta}:=U_{\delta}gf_{\varepsilon}$,
for every $g \in SO(3)$ and $\varepsilon>0.$ We denote this
last measure by $\nu_{g,\eps,\delta}$.  Let
$R(g,\eps,\delta) = R(\nu_{g,\eps,\delta})$.

\begin{definition} The {\em ($\delta$-)diffused random exponent}
is the average over $SO(3)$
of $R(g,\eps,\delta)$:
$$ R(\eps,\delta)=\int_{g \in SO(3)}R({g,\eps,\delta})\dd \nu(g) .$$
\end{definition}

\begin{lemma}\label{l=limit}
$$\limsup_{\delta\rightarrow 0}\int_{g \in SO(3)}R({g,\eps,\delta})\dd g
\leq \Lambda(\eps).$$
\end{lemma}
\begin{proof}
Note that $\lim_{\delta\to 0} \nu_{g,\eps,\delta} = \delta_{gf_\eps}$,
Dirac
measure supported on $gf_\eps$.
By the previous lemma,
\begin{eqnarray*}
\limsup_{\delta\rightarrow
0}R({g,\eps,\delta}) &=&
\limsup_{\delta\rightarrow
0}R(\nu_{g,\eps,\delta}) \\
&\leq& R(\delta_{gf_\eps} )\\
&=& \lambda(gf_{\varepsilon})
\end{eqnarray*}
for each $g \in SO(3)$, so the same is true for the integral.
\eproof\end{proof}

Now let $h_0:\Sset^2\to \Sset^2$ be any $\mu$-preserving diffeomorphism,
and let $\H=U_{\delta}h_0$.  As in the previous section,
define the space
$$\H_\sharp = \{h_\sharp\,\mid\, h\in \H\},$$
and evaluation maps
$$ev: \H\times \Sset^2\to \Sset^2,\qquad ev_\sharp: \H\times P\Sset^2\to
P\Sset^2.$$
Let $\nu_\delta$ and $\nu_{\delta\sharp}$ be the push-forwards of
normalized
Haar measure to $\H$ and $\H_\sharp$, respectively.
Note that $\mu$ is stationary for the process generated by
$(\H, \nu_\delta)$.

\begin{prop}\label{p=formula} If $h_0$ is not an isometry, then
for fixed $\delta>0$
the random process on $P\Sset^2$
generated by $(\H_\sharp, \nu_{\delta\sharp})$
has a stationary measure $m_{\delta}$
that is absolutely continuous with smooth density, covers $\mu$,
and is the unique such stationary measure.

Moreover,
$$R(\nu_{\delta})=\int_{P\Sset^2}\log\|Th_0 v\|\dd m_{\delta}(v)>0$$
\end{prop}

\begin{cor}\label{c=formula}
For fixed $\delta>0$, $\eps\neq 0$ and $g\in SO(3)$,
the random process on $P\Sset^2$
generated by $(\F_{g,\eps,\delta},\nu_{g,\eps,\delta})$
has a stationary measure $m_{g,\eps,\delta}$
that is absolutely continuous with smooth density, covers $\mu$,
and is the unique such stationary measure.

Moreover,
$$R(g,\eps,\delta)=\int_{P\Sset^2}\log\|Tf_\eps v\|\dd m_{g,\eps,\delta}>0$$
\end{cor}

\noindent{\bf Proof of Corollary~\ref{c=formula}:}\,
Apply Proposition~\ref{p=formula} to the case where
$h_0=gf_\eps$.

\bigskip

\noindent{\bf Proof of Proposition~\ref{p=formula}:}\,
Let $h_0$ and $\delta>0$ be given.  We break the proof into steps.

\bigskip

\noindent{\bf Step 1: Construction of $m_{\delta}$.}

Recall that the convolution of a probability measure $\pi$ on
$\H_\sharp$ and a probability
measure $m$ on $P\Sset^2$ is a probability measure $\pi\star m$ on
$P\Sset^2$ defined by
$$\pi \star m(E) = \int_{\H_\sharp}m(h_\sharp^{-1}E)\dd\pi(h_\sharp)$$ for
every $m$-measurable $E\subseteq P\Sset^2$.
That a measure $m$ is stationary for the measure
$\pi$ is equivalent to the fact that $\pi_\sharp \star m = m$.
For $k>1$ we let $\pi^k \star m = \pi \star(\pi^{k-1} \star m)$.
For any probability measure $m$ on $P\Sset^2$ any weak limit of the
Ces\`aro sums $\frac{1}{n}\sum _1^n \pi^k \star m$
is a stationary measure for $\pi$.
Beginning with a measure $m$ which pushes forward under projection to
$\mu$ produces
an invariant measure by this process with the same property. If we start
with $m$
as Liouville measure  on $P\Sset^2$ and
$\pi = \nu_{\delta\sharp}$ we call this limiting measure $m_{\delta}$.

\bigskip

\noindent{\bf Step 2:  $m_{\delta}$ is absolutely continuous, with smooth
density}

For any measurable set $A\subseteq P\Sset^2$, we have:
\begin{eqnarray*}
m_{\delta}(A)&=&(\nu_{\delta\sharp}\times m_{\delta}) ev_\sharp^{-1}(A)\\
&=&(\nu_{\delta\sharp}\times m_{\delta}) \{(h_\sharp,v)|h_\sharp(v) \in
A\}\\
&=&\int_{v \in P\Sset^2}\nu_{\delta\sharp}\{h_\sharp\,|\,h_\sharp(v) \in
A\}\dd m_{\delta}(v)
\end{eqnarray*}

Now if the Liouville measure $m(A)$ equals zero, then
$\nu_{\delta}\{h_\sharp|h_\sharp(v) \in A\}$ must also be zero,
for every  $v \in P\Sset^2$.  Thus  $m_{\delta}(A)$ is zero.
It follows that $m_{\delta}(A)$ is absolutely
continuous with respect to $m$. So there is a non-negative integrable
function
$\varphi_{\delta}$ defined on $P\Sset^2$ so that for any measurable  $A
\subseteq P\Sset^2$,
$$m_{\delta}(A)= \int _A
\varphi_{\delta}(x) dm(x).$$
In Lemma \ref{l=density} we will prove that $\varphi_{\delta}$
satisfies the following formula:
\begin{eqnarray*}
\varphi_{\delta}(x)&=&\frac{1}{m(B(x,\delta))}\int_{h_{0\sharp}^{-1}B(x,\delta)}\varphi_{\delta}(
y)\dd m(y)\\
&=&\frac{1}{m(B(x,\delta))}\int_{B(x,\delta)}\varphi_{\delta}(
h_{0\sharp} ^{-1}z ) \hbox{Jac}(h_{0\sharp}^{-1}) (z)\dd m(z).
\end{eqnarray*}
It follows now fairly directly that $\varphi_{\delta}$ is  as smooth
as $h_{0}$, since the average over a $\delta$-ball of an $L^1$ function
is continuous, and of a $C^k$ function, is $C^{k+1}$.

\bigskip

\noindent{\bf Step 3: $R(\nu_\delta)$ satisfies the integral formula,
and the exponents of $m_{\delta}$ are nonzero.}

The argument that
$R({\nu_{\delta}})=\int_{P\Sset^2}\log\|Th_0 v\|\dd m_{\delta}$
is now the same as in the proof of proposition \ref{P=random}, using
Birkhoff's
and Oseledec's theorems, where
$\mathcal F_\sharp$ is replaced by $\H_\sharp$, $\nu$ by $\nu_{\delta}$
and $m$ by
$m_{\delta}$.

Next we will prove that the largest exponent is positive and from
that we will deduce uniqueness. As in Proposition~\ref{P=random}, for
any $h = gh_0 \in\H$ the Jacobian
of $h_\sharp$ with respect to $m_{\delta}$  at the vector $v$ in
$T_{1}\Sset^2$ is
\begin{eqnarray}\label{e=density}
\rho(h,v) &=&
\frac{\varphi_{\delta}(h_{\sharp} v)}{\|Th v\|^2 \varphi_{\delta}(v)}
 = \frac{\varphi_{\delta}(gh_{0\sharp} v)}{\|Th_0 v\|^2
\varphi_{\delta}(v)}
\end{eqnarray}
provided that $\varphi_\delta(v)\neq 0$.

We claim that the function
$\rho(h,v)$ cannot be $\nu_\delta\times m_\delta$ - almost everywhere
equal to $1$.  Suppose for the sake of contradiction
that $\rho(h,v)=1$ a.e. Since $\varphi_\delta$ is continuous,
if we fix $m_\delta$-a.e. $v$, then the function $\rho(\cdot, v)$ is
continuous. Thus, for $m_\delta$-a.e. $v$, we must have $\rho(h,v)=1$ for
{\em every} $h\in U_\delta$.

 Next, notice in the expression (\ref{e=density}) for $\rho(h,v)$,
that the only term that depends on $g\in U_\delta$ is the numerator
$\varphi_{\delta}(gh_{0\sharp} v)$.  Rewriting this
expression, we have, for almost every  $v$
in the support of $\varphi_\delta$,
\begin{eqnarray}\label{e=const}
{\varphi_{\delta}(gh_{0\sharp} v)} &=&  \|Th_0 v\|^2 \varphi_{\delta}(v),
\end{eqnarray}
for every $g\in U_\delta$.
Since
$m_\delta$ projects to $\mu$, and
$\varphi_\delta$ is continuous,
we have that for $\mu$-a.e. $x\in \Sset^2$, the set
$O =  \{v\,\mid\, \varphi_\delta(v) > 0\}$ is an open,
$h_{0\sharp}$-invariant set in $P\Sset^2$
that intersects almost every fiber.
Equation (\ref{e=const}) implies
$\varphi_\delta$ must be constant on connected components of $O$, since
varying $g$ inside of $U_\delta$, the vector $gh_{0\sharp} v$ covers
an open neighborhood of $h_{0\sharp} v$  in $P\Sset^2$.

But, again by equation (\ref{e=const}), on each such component of $O$, we
must have that
$\|Th_0\|$ is constant.  Since $O$ intersects almost every fiber of
$P\Sset^2$, we obtain that for almost every $x\in \Sset^2$, there exists
a connected open set $I_x$ in the fiber of $P\Sset^2$ on
which $T_xh_0$ has constant norm.  Since $h_0$ preserves area,
we must have $\|T_xh_0\| = 1$ for $\mu$-a.e. $x$, contradicting the
assumption that $h_0$ is not an isometry.

So $\rho(h,v)$ is not a.e. equal to $1$,  and by Jensen's inequality,
we have:
\begin{eqnarray*}
\int_{\H\times P\Sset^2} \log \rho(h,v) \dd\nu_\delta(h)\,dm_\delta(v)
& < & \log\int_{\H\times P\Sset^2} \rho(h,v)
\dd\nu_\delta(h)\,dm_\delta(v)\\
&=&0.
\end{eqnarray*}

But
\begin{eqnarray*}
R(\nu_\delta) & = & \int_{P\Sset^2}\log\|Th_0 v\|\dd m_\delta(v)\\
& = & -\frac{1}{2} \int_{\H\times P\Sset^2}\log\|Th v\|^{-2}\dd
\nu_\delta(h)\,dm_\delta(v)\\
& =&  -\frac{1}{2} \int_{\H\times P\Sset^2}\log\left(\|Th
v\|^{-2}\frac{\varphi_\delta(h_\sharp(v))}{\varphi_\delta(v)}\right)\dd\nu_\delta(
h)\dd m_\delta(v),\\
&=& -\frac{1}{2}\int_{\H\times P\Sset^2} \log \rho(h,v)
\dd\nu_\delta(h)\dd m_\delta(v)\\
&>& 0.
\end{eqnarray*}
(Here we used the stationarity of $m_\delta$ to conclude that
the integral of $\log({\varphi_\delta(h_\sharp(v))}/{\varphi_\delta(v)})$ is
$0$).

\bigskip

\noindent{\bf Part 4:  $m_{\delta}(A)$  is unique}

It remains to prove that the measure $m_{\delta}$ is unique among
absolutely continuous stationary
measures which cover $\mu$, now that we know that the
random Lyapunov exponents are not zero.

Let $\gamma$ be any such
measure. Then, as for  $m_{\delta}$,
there is a nonnegative smooth function $\psi$ such that for
any measurable $A \subseteq P\Sset^2$,
$$\gamma(A) = \int_A\psi(v)\dd m(v).$$
Let $\gamma_x$ be the disintegration of $\gamma$ on the fiber
$T_{1,x}\Sset^2$.  The density of $\gamma_x$ with respect to
Lebesgue measure on the fiber $T_{1,x}\Sset^2$ is the restriction of
$\psi$ to the fiber.

We need the notion of a natural extension of a non-invertible
transformation.
Let $(\Omega, {\mathcal A}, p)$ be a
probability space with $p$-preserving transformation
$T:\Omega\to \Omega$.  Let
$$\hat\Omega = \{(\cdots, \omega_{-1}, \omega_0)\in
\Pi_{-\infty}^{1}\Omega\,\mid\, T(\omega_{-i})
= \omega_{-i+1},\,\forall i\geq 1\}.$$
Let $\hat T: \hat\Omega\hookleftarrow$
be the map:
$$\hat T(\cdots, \omega_{-1}, \omega_0)
 = (\cdots, \omega_{-1}, \omega_0, T(\omega_0)),$$
and let $\hat\pi_0: \hat{\Omega}\to\Omega$
be the projection onto the first factor:
$$\hat\pi_0(\cdots, \omega_{-1}, \omega_0) = \omega_0.$$
Let $\hat{\mathcal A}$ be the smallest $\sigma$-algebra
on $\hat\Omega$ so that $\hat\pi_0$ and $\hat T $ are both measurable.
On $(\hat\Omega,\hat{\mathcal A})$, there is a unique probability
measure $\hat p$, invariant under $\hat T $, that
pushes forward under $\hat\pi$ to $p$.  The measure-preserving
system $(\hat\Omega,\hat{\mathcal A},\hat p)\hookleftarrow \hat T $ is
called
the {\em natural extension} of $(\Omega,{\mathcal A},p)\hookleftarrow  T $.
The natural extension $\hat T $ is invertible, and ergodic
if $ T $ is ergodic.

Let $\tau$ and $\tau_\sharp$ be the associated
auxiliary maps to the processes generated by
$(\H,\nu_\delta)$ and $(\H_\sharp,\nu_{\delta\sharp})$,
respectively.
The natural extension
of $\tau_\sharp$ with respect to the measure
$\nu_{\sharp\delta}^\infty \times \gamma$ fibers over
the natural extension of $\tau$ with respect to $\nu_{\delta}^\infty\times
\mu$; the fiber over $((h_{i})_{-\infty}^\infty,x)$  is $T_{1,x}\Sset^2$.
For $\hat\nu^\infty$-almost every ${\bf h}=(h_{i})_{-\infty}^\infty$
and almost every fiber $T_{1,x}\Sset^2$ there is a measure
$\hat\gamma_{{\bf h}, x}$
which is the disintegration of
$\widehat{\nu^\infty \times \gamma}$ along the fiber.
Note that the extension $\widehat{\nu^\infty\times \gamma}$ is
determined by this system of measures, and therefore so is $\gamma$.

By (\cite{L-Q}, Proposition 1.1, p.131), if
$({\bf h}=(h_{i})_{-\infty}^\infty ,x) \in \H_{-\infty}^\infty \times
\Sset^2 =
\widehat{\H^\infty \times \Sset^2 } $,
then $\hat\gamma_x$ is the limit of the push-forwards:
\begin{eqnarray}\label{e=limit}
\hat\gamma_{{\bf h}, x}& = &
\lim_{n\to\infty}(h_{-n\sharp}\dots h
_{1\sharp})_\ast(\gamma_{h_{-n}^{-1} \dots h_1^{-1}(x)})
\end{eqnarray}

\begin{lemma}
The limit (\ref{e=limit}) does not depend on the initial choice
of absolutely continuous stationary measure $\gamma$ covering $\mu$.
\end{lemma}
\begin{proof}
The average exponents of $\tau$ are nonzero, so the average exponents of
$\hat\tau$ are also nonzero.
It is easy to see that $\tau$ is ergodic,
and thus, so is $\hat\tau$ and the exponents of $\hat\tau$ are
in fact nonzero $\mu$-a.e.

Let $u({\bf h}, x)\in P\Sset^2$ be the unstable Lyapunov direction
for $\hat\tau$ over $({\bf h}, x)$. For every
$\epsilon>0$, there exists an $n>0$, and a set $G\subset
\H_{-\infty}^\infty \times \Sset^2$ such that
\begin{itemize}
\item $\nu_{-\infty}^\infty\mu(G) > 1-\epsilon$
\item for every $({\bf h},x)\in G$,
the $\gamma_{{\bf h},x}$-measure of
an $\epsilon$-neighborhood of  $u({\bf h}, x)$ in
$T_{1,x}\Sset^2$ is at least $1-\epsilon$.
\end{itemize}
It follows that the limit in (\ref{e=limit}) is concentrated
on the point  $u({\bf h}, x)$.
\eproof\end{proof}

Thus the natural extension of $\tau$ and $\gamma$ is the same as the
natural extension of $\tau$ and  $m_{\delta}$,
so $\tau$  and $m_{\delta}$ are themselves equal.
\eproof

\bigskip

The next lemma completes the proof of Proposition~\ref{p=formula}.

\begin{lemma} \label{l=density} Let $\gamma$ be an absolutely continuous
measure with respect to $m$ on $P\Sset^2$
that is stationary for $\nu_{\delta}$.
Then the density function $\psi$ defining $\gamma$ satisfies
$$\psi(x)=
\frac{1}{m(B(x,\delta))}\int_{h_{0\sharp}^{-1}B(x,\delta)}\psi(y)\dd
m(y).$$
\end{lemma}

\begin{proof}
For any measurable set $A \subseteq P\Sset^2$, we have:
\begin{eqnarray*}
\int_A\psi(y)\dd m(y)&=&\gamma(A)\\
&=&\nu_{\delta} \times \gamma (ev_\sharp^{-1}(A))\\
&=&\nu_{\delta} \times \gamma \{(h,y) \in \H_{\sharp}\times
P\Sset^2\,\mid\,h(y) \in A\}\\
&=&\int_{g \in U_\delta }\gamma\{h_{0\sharp}^{-1}g^{-1}A\}
\dd\nu_{\delta}(g)\\
\end{eqnarray*}
\begin{eqnarray*}
&= &\int_{g \in U_\delta }\int_{h_{0\sharp}^{-1}g^{-1}A}\psi(y)\dd
m(y)\dd\nu_{\delta}(g)\\
&= & \int_{g \in U_\delta
}\int_A\psi(h_{0\sharp}^{-1}g^{-1}(x))J(h_{0\sharp}^{-1}g^{-1}(x))^{-1}
\dd m(x)\dd\nu_{\delta}(g)\\
&=&\int_A \int_{g \in U_\delta
}\psi(h_{0\sharp}^{-1}g^{-1}(x))J(h_{0\sharp}^{-1}g^{-1}(x))^{-1}
\dd\nu_{\delta}(g)\dd m(x)
\end{eqnarray*}
So,
$$\psi(x)=\int_{g \in U_\delta
}\psi(h_{0\sharp}^{-1}g^{-1}(x)J(h_{0\sharp}^{-1}g^{-1}(x))^{-1}
\dd\nu_{\delta}(g).$$
Since the integrand only depends on the point $x\in P\Sset^2$ we push
the measure
$\nu_{\delta}$ forward to
$P\Sset^2$, use that $J(g)^{-1} = 1$ and that Haar measure on $SO(3)$
pushes forward to Liouville measure on $P\Sset^2$
to obtain
$$\psi(x) = \frac{1}{m(B(x,\delta))}\int_{y\in B(x,\delta)}
\psi(h_{0\sharp}^{-1}(y))J(h_{0\sharp}^{-1}(y))\dd m(y).$$
Finally, changing variables one more time, gives:
$$\psi(x) =
\frac{1}{m(B(x,\delta))}\int_{h_{0\sharp}^{-1}B(x,\delta)}\psi(y)\dd
m(y).$$
\!\!\eproof

This completes the proof of Proposition~\ref{p=formula}.

\end{proof}

\bigskip

\noindent

Returning  to discussion of the family $\F_\eps$, we have
verified that properties 1.-3. hold for the measures $m_{g,\eps,\delta}$.
Now let
$$m_{\eps,\delta}=\int_{g \in
SO(3)}m_{g,\eps,\delta}\,d g.$$
Since
$$\int_{SO(3)} R(g,\eps,\delta)\, dg = \int_{P\Sset^2} \log\|Tf v\| \,dm_{g,\eps,\delta}(v)\,
dg,$$
it follows from Corollary~\ref{c=formula} and Lemma~\ref{l=limit} that:

\begin{prop}\label{P=theinequality}
$$\Lambda(\eps)\geq
\limsup_{\delta\to0}\int_{P\Sset^2}\log\|Tf_\varepsilon(x)v\|\dd
m_{\eps,\delta}$$
\end{prop}

\subsection{What is this measure $m_{\epsilon,\delta}$?}

For the family of twist maps $\F_\eps$ under consideration, we prove
in  Section~\ref{S=smalleps} that the answer to Question~\ref{Q6} is ``no'',
and for small $\eps$, we have $\Lambda(\eps) <  R(\eps)$.  It then follows from
Proposition~\ref{P=theinequality} that
for small $\eps$,  Lebesgue measure $m$ is not a weak limit of
$m_{\eps,\delta}$ as $\delta\to 0$.  At the opposite extreme,
we show in Section~\ref{S=epstoinfty} that
as $\eps$ tends to infinity,
the measures $m_{\eps,\delta}$ do approach Lebesgue measure, for
$\delta>0$ fixed.
We hope that a future experiment will reveal more precisely
how these measures behave in $\delta$, for moderate values of $\eps$.
As will be seen in section \ref{S=Numerics}, the tiny differences
between $\Lambda(\eps)$ and $R(\eps)$ for moderate and large values of
$\eps$, are expected to give rise to numerical difficulties in estimating
the behavior with respect to $\delta.$

In Section~\ref{S=linear}, we show that in $SO(2)$-invariant families
of $2\times 2$ matrices, if $m_\delta$ is the analogous ``in-between''
measure, averaged over the family, then $m_\delta$ is Lebesgue measure
on $S^1$, for all $\delta>0$.   On the other hand, we also
show that for unitarily-invariant
families, as $\delta\to 0$,
the $m_\delta$  do {\em not} approach the natural unitarily invariant
measure on the appropriate Grassmannian manifold,
but instead they  limit on an even ``better'' measure, in the sense that
this measure forces
the (strict) inequality $\Lambda > R$.   We describe the construction of
$m_\delta$ for matrices in Section~\ref{S=linear}.

\section{The linear case}\label{S=linear}
Question~\ref{Q6} was originally motivated by a result of Dedieu-Shub about
random and deterministic exponents for families of matrices.
In this section, we describe these results and apply the framework of
the previous section to the matrix setting.

Let $L_i$ be a sequence of linear maps  mapping finite
dimensional normed vector spaces $V_i$ to $V_{i+1}$ for $i \in
\bf{N}$. Let $v \in V_0\backslash\{0\}$. If the limit $\lim
\frac{1}{k} \log \|  L_{k-1} \ldots L_0(v) \|$ exists it is
called a Lyapunov exponent of the sequence. It is easy to see
that if two vectors have the same exponent then so does every
vector in the space spanned by them. It follows that there are at
most $\dim(V_0)$ exponents.  We denote them $\lambda_j$  where $j
\leq k \leq \dim(V_0)$. We order the $\lambda_i$ so that
$\lambda_i \geq \lambda_{i+1}$.

Given a probability measure $\mu$ on $\GL$,
the space of invertible $n \times n$
complex matrices, we may form infinite sequences of elements
chosen at random from $\mu$ by taking the product measure on
$\GL^{\bf N}$. Thus we may
also talk about the Lyapunov exponents of sequences or almost all
sequences in $\GL^{\bf N}$.

For measures $\mu$ on $\GL$
satisfying a mild integrability condition, we have, by
Oseledec's Theorem, $n$ Lyapunov
exponents  $r_1 \geq r_2 \geq \ldots \geq r_n \geq - \infty$ such
that, for almost every sequence $\ldots g_k \ldots g_1 \in \GL$,
the limit
$$\lim \frac{1}{k} \log \|  g_k \ldots g_1v \|$$ exists
for every $v \in \Cn \setminus \{ 0 \}$ and equals one of the
$r_i$, $i = 1 \ldots n$, see Gol'dsheid and Margulis \cite{G-M}
or Ruelle \cite{R} or Oseledec \cite{O}. We may call the numbers
$r_1,\ldots , r_n$ random Lyapunov exponents or even just random
exponents. If the measure is concentrated on a point $A$, these
numbers:
$$\lambda_i(A) = \lim \frac{1}{n}\log \|  A^n v \|, \qquad i=1\ldots n,$$
are $\log \vert
e_1 \vert$, \ldots , $\log \vert e_n \vert$, where
$e_i (A) = e_i$, $i = 1 \ldots n$, are the
eigenvalues of $A$ written with multiplicity and $\vert e_1
\vert \geq \vert e_2 \vert \geq \ldots \geq \vert e_n
\vert .$

The integrability condition for Oseledec's Theorem is
$$g \in \GL \rightarrow \log^+( \| g \| )\ \mbox{is }\mu-\mbox{integrable}$$
where for a real valued function $f$, $f^+ = \max [0,f].$ Here we
will assume more so that all our integrals are defined and
finite, namely:
$$(*)\hskip 2cm g \in \GL \rightarrow \log^+( \| g \| )\ \mbox{and }
\log^+( \| g^{-1} \| )\ \mbox{are }\mu-\mbox{integrable}.\hskip 2cm$$

In this matrix setting, there are rigorous lower bounds
for the average exponents (= logarithms of moduli of eigenvalues)
of unitarily-invariant families in $\GL$.  In  \cite{DeSh},
the following bound is proved:

\begin{theorem}{\rm \cite{DeSh}}
\label{th1} If $\mu$ is a unitarily invariant measure on $\GL$
satisfying $(*)$ then, for $k=1, \ldots , n,$
$$\int_{A \in \GL} \displaystyle \sum_{i=1}^k \log \vert \lambda_i (A)
\vert d\mu (A) \geq \displaystyle \sum_{i=1}^k r_i.$$
\end{theorem}

By unitary invariance we mean $\mu (U(X)) = \mu (X)$ for all
unitary transformations $U \in \Un$ and all $\mu$-measurable $X
\subset \GL$.

We can rephrase a special case of this
theorem in a form similar to Question~\ref{Q6}.
Fix $A\in\GL$.
As above, let $\nu$ be normalized Haar measure on $\Un$,
and also denote by $\nu$ the push-forward of $\nu$
to the coset $\Un A\subset \GL$. Let $R(A) = r_1(A)$
be the largest random exponent of $\nu$, and let
$$\Lambda(A) = \int_{B\in \Un A} \log|e_1(B)|\dd \nu(B).$$
Then we have:

\begin{cor}{\rm \cite{DeSh}}\label{C=desh} For $n\geq 2$, and for any $A\in\GL$,
$$\Lambda(A) \geq R(A).$$
Equality holds if and only if $A\in \Un$.
\end{cor}

Thus non-zero Lyapunov exponents for the family, i.e. non-zero
random exponents, implies that at least some of the individual
linear maps have non-zero exponents,  i.e eigenvalues of modulus not
equal to 1. Hence the question that we posed for diffeomorphisms
has a positive answer for sufficiently rich (i.e. unitarily-invariant)
families of matrices.

\begin{remark} Theorem \ref{th1} is not true for general measures on $\GL$ or
$\GLR$ even for $n=2$. Consider
$$ A_1 = \left ( \displaystyle
\begin{array}{cc}
1&0\\
1&1\\
\end{array} \right ), \ \
A_2 = \left ( \displaystyle
\begin{array}{cc}
1&1\\
0&1\\
\end{array} \right ),$$
and give probability $1/2$ to each. The left hand integral is zero
but as is easily seen the right hand sum is positive. So, in this
case the inequality goes the other way. We do not know a
characterization of measures which make Theorem \ref{th1} valid.
\end{remark}

We expect similar results for orthogonally invariant probability
measures on $\GLR$ but the only case in which such a result has
been proved is in dimension
2, where we have:

\begin{theorem} \label{th3}{\rm \cite{DeSh}}
Let $\mu$ be a probability measure on $\GLDR$ satisfying
$$\hskip 3cm g \in \GLDR \rightarrow \log^+( \| g \| )\quad \mbox{and}\quad
\log^+( \| g^{-1} \| )\quad \mbox{are }\mu-\mbox{integrable}.\hskip
5mm$$

\noindent a. If $\mu$ is a $\SO2$ invariant measure on $\GLDPR$
then,
$$\int_{A \in \GLDPR} \displaystyle \log \vert \lambda_1 (A)
\vert d\mu (A) = \int_{A \in \GLDPR} \int_{x \in \S1} \log \| Ax
\| \dd x \dd \mu (A).$$

\noindent b. If $\mu$ is a $\SO2$ invariant measure on $\GLDMR$,
whose support is not contained in $\RE \O2$ i.e. in the set of
scalar multiples of orthogonal matrices, then
$$\int_{A \in \GLDMR} \displaystyle \log \vert \lambda_1 (A)
\vert d\mu (A) > \int_{A \in \GLDMR} \int_{x \in \S1} \log \| Ax
\|  \dd x \dd\mu (A).$$
\end{theorem}

Here $\GLDPR$ (resp. $\GLDMR$) is the set of invertible matrices
with positive (resp. negative) determinant.
To rephrase Theorem 3 a., fix $A\in \GLDPR$, and let
$\nu$ be the push forward of Haar measure on $\SO2$ to the coset
$\SO2 A$.  Let $R(A)$ be the largest random exponent for
the process induced by $\nu$, and let
$\Lambda(A) = \int_{B\in \SO2 A} \log|e_1(A)|\dd\nu(A)$.
The we have:

\begin{cor}{\rm \cite{DeSh, AvBo}}\label{C=so2} For any $A\in\GLDPR$,
$$\Lambda(A) = R(A).$$
\end{cor}

We give an alternate proof of Corollary~\ref{C=so2} in the following
subsection.

\subsection{$m_\delta$ for matrices}

Let $A \in \GL $ or $\GLR$  and $\mu$ be the Haar measure on $\Un$
or $SO(n,{\bf R})$, respectively, normalized to be a probability measure.
Let $G$ denote $\GL $ or $\GLR$.  As we did for families
of diffeomorphisms in Section 3,  we now interpolate between random
products and deterministic powers of matrices by changing $\mu$. Let
$\delta>0$ and  $G_\delta$ be the $\delta$-neighborhood of the
identity in $G$. For $g \in G$, $G_\delta g A$ is a neighborhood
of $gA$ in $GA$. We normalize Haar measure restricted to
$G_\delta$ and push it forward to $G_\delta g A$. Let us call this
measure $\mu_{\delta,g}.$ Let $r_1(\delta,g)$ be the largest
random exponent for this measure.  At the end of this subsection, we prove:

\begin{prop}\label{P=limitmeasures}
$\lim_{\delta\rightarrow 0}r_1(\delta,g)=\log|e_1(gA)|.$
\end{prop}

Let $\GR$ denote the Grassmannian manifold of $k$ dimensional
vector subspaces in $\C^n$, and let $m$  be the natural unitarily
invariant probability measure on $\GR$. Any $n\times n$ complex
matrix acts on the homogeneous space $G_{n,k}$ by
left-multiplication. For $A \in \GL$  and $P\in\GR$,  denote by
$A|P$ the restriction of $A$ to the subspace $P$. Now let
$m_{\delta,g}$ be the stationary measure on the Grassmannian
$G_{n,1}(={\bf C}P^{n-1})$ induced by $\mu_{\delta,g}.$

\begin{prop}\label{P=randomeponent}
$r_1(\delta,g)=\int_{P \in G_{n,1}} \log  \|
(A|P)\|\dd m_{\delta,g}$
\end{prop}

Let $m_\delta=\int_{g \in G}\nu_{\delta,g}d \mu.$ It follows that:

\begin{prop}\label{P=limitformula}
 $\int_{g \in G}\lambda_1(gA)d \mu=\lim_{\delta\rightarrow
0}\int_{P \in G_{n,1}}\log \| (A|P)\|\dd m_{\delta}.$
\end{prop}

Now
$$r_1=\int_{P \in G_{n,1}}\log \| (A|P)\|\dd m(P).$$
So a comparison of $\int_{g \in G}\lambda_1(gA)d \mu$
and $r_1$ can be achieved via an understanding of the relationship
between $m$ and $m_{\delta}.$
We have two results in this direction.

First, recall that inequality in Corollary~\ref{C=desh}
is strict when $n \geq 2$ unless $A$ is an isometry.
So for
$G=SU(n)$ the measures $m_{\delta}$ favor the expanding
directions of $A$ as $\delta\rightarrow0$.
By Proposition~\ref{P=limitformula}, we obtain the immediate
corollary:

\begin{cor}
For $G=SU(n)$, $n\geq 2$,
$$\lim_{\delta\to 0} m_{\delta} \neq m,$$
unless $A$ is an isometry.
\end{cor}

Experimentally the same
seems to hold for $SO(n)$ when $n>2$, but we
have not checked this very carefully.

By contrast, the equality in Corollary~\ref{C=so2} is consistent with
$\lim_{\delta\to 0} m_{\delta} = m$. In fact, more is
true:

\begin{theorem}\label{T=equalmeasures}
For $G=SO(2)$, $$m_{\delta}=m,$$ for all
$\delta>0.$
\end{theorem}

Combined with Proposition~\ref{P=limitformula},
Theorem~\ref{T=equalmeasures} gives
another proof of Corollary~\ref{C=so2}.

\bigskip

\noindent{\bf Proof of Proposition~\ref{P=limitmeasures}:}
Following the proof of  Lemma~\ref{l=limit}, one obtains that
$$\log|e_1(gA)| \geq \lim_{\delta\to 0}  r_1(\delta, g).$$

If $A$ is replaced by $cA$, for $c\in\C\setminus\{0\}$,
then both sides of the equality change by $\log|c|$.
So we may assume that $|\det A| = 1$, and it will be enough to prove that
$\log|e_1(gA)| \leq \lim_{\delta\to 0} r_1(\delta, g)$
under the hypothesis that $|e_1(gA)|>1$.

Let $E\subset\Cn$ be the generalized eigenspace of the eigenvalues of
$gA$ whose modulus equals $|e_1(gA)|$.  Then given $\eps>0$, there
is a metric on $\Cn$, a closed cone $K\subset \Cn$ containing $E$ in
its interior, and a $\delta>0$ such that, for any $B\in\GL$ in
the $\delta$-neighborhood of $gA$, we have:
\begin{enumerate}
\item $B(K)\subset K$, and
\item $\|Bv\| \geq (|e_1(gA)|-\eps) \|v\|$, for all $v\in K$.
\end{enumerate}
It follows that, for any $\eps>0$, there is a $\delta>0$ such that,
$$\|B_n\cdots B_1\| \geq (|e_1(gA)|-\eps)^n,$$
for all sequences $B_1,\ldots, B_n$ in $U_\delta g A$.  Hence
$$r_1 (\delta, g)  \geq \log(|e_1(gA)|-\eps),$$
and so
$$\lim_{\delta\to 0} r_1(\delta, g) \geq \log|e_1(ga)|.$$
\!\!\eproof

\bigskip

We now turn to the proof of Theorem~\ref{T=equalmeasures}.

\bigskip

\noindent{\bf Proof of Theorem~\ref{T=equalmeasures}:}
By an argument presented in the proof of Proposition~\ref{P=limitmeasures},
we may assume that $\det A = 1$.  The projective action of
$SL(2,{\bf R})$ on ${\bf R}P^1$ is conjugate to the standard action
on the circle $S^1 = \{z\in \C\,\vert\, |z|=1\}$
by linear fractional transformations.
The conjugacy sends the rotation by $\theta$ to rotation by $2\theta$.
Let $f:S^1\to S^1$ be the linear fractional transformation
induced by $A$, and let $\F = \{\alpha f\,\vert\, \alpha\in S^1\}$.
Let
$$U_{\alpha,\delta} = \{\beta \alpha f\,\vert\, \arg(\beta)\in (-\delta,\delta)\},$$
and let $\nu_{\alpha,\delta}$ be normalized
Lebesgue measure on $U_{\alpha,\delta}$,
pushed forward from $(-\delta,\delta)$.
Denote by $m_{\alpha,\delta}$ the stationary measure on
$S^1$ induced by $\nu_{\alpha,\delta}$.
We will show that:
$$m_\delta = \int_{S^1} m_{\alpha,\delta} \dd \alpha$$
is Lebesgue measure on $S^1$.

The same argument as in the proof of Lemma~\ref{l=density} shows
that for every $\alpha\in S^1$,
$$dm_{\alpha,\delta}(z) = \varphi_{\alpha,\delta}(z) dz,$$
where
 $$\varphi_{\alpha,\delta}(z)=
\frac{1}{2\delta}\int_{y\in (\alpha f)^{-1}B(z,\delta)}\varphi_{\alpha,\delta}(
y)\dd y,$$
and
$B(z,\delta) = \{\beta z\,\vert\, \arg(\beta)\in (-\delta,\delta)\}$.
Setting
$$k_\delta(z) = \frac{1}{2\delta} 1_{B(1,\delta)}(z),$$
we have that
$$\varphi_{\alpha,\delta}(z)=
\int_{y\in S^1} k_\delta(\alpha \overline{z} f(y)) \varphi_{\alpha,\delta}(
y)\dd y,$$
where we use $\overline{z}$ to denote the multiplicative inverse of
$z\in S^1$. Note that $\int_{S^1} k_\delta(z) \dd z = 1$.

Consider the following, more general setting.
Let $k $ be a non-negative function on $S^1$ such that $\int k = 1 $ and $\int
k^2 <
\infty $ (all integrals are with  respect to the normalized Lebesgue measure on
$S^1$). For example, $k=k_\delta$.
Let $f: S^1 \mapsto S^1$ be a fractional linear transformation, so that we can
write, for $\vert z \vert = 1$: $$ f(z) = \sum _{n\geq 0} c_n z^n .$$

For $\alpha \in S^1 $ , define the operator $ L_{\alpha, f} $ on
real functions on $S^1$
by: $$ L_{\alpha, f} \varphi (z) = \int k(\alpha \overline{z}
{\mathstrut f (y)} )
\varphi (y) dy.$$

The operator $L_{\alpha, f}$ is a positive operator, $\int L_{\alpha, f}
\varphi = \int \varphi $. There exists a unique function $\varphi _\alpha $
satisfying $L_{\alpha, f} \varphi _\alpha = \varphi _\alpha $ and $\int \varphi
_ \alpha =1$. The function $\varphi _\alpha $ is positive and continuous,
upper and positive lower bounds for $\varphi _ \alpha $ can be chosen
uniformly in  $\alpha $.  In the case $k = k_{\delta}$, we have
$\varphi_\alpha = \varphi_{\alpha,\delta}$.

\begin{prop}
We have, for all $z \in S^1$,  $\int \varphi _\alpha (z) d\alpha =  1$.
\end{prop}

The proposition follows directly from the  following two claims:

\bigskip

\noindent{\bf Claim 1.}  For all $m\geq 0$, all $z \in S^1$,  $\int L_{\alpha, f}^m 1
(z) d\alpha =1 $.

\bigskip

\noindent{\bf Claim 2.} The sequence $\frac{1}{n} \sum _{m = 1}^n L_{\alpha, f}^m 1
(z)$ converges to $\varphi _\alpha (z)$ in $L^1(d\alpha, dz)$ as $n\to
\infty $.

\bigskip

Claim 2 follows from the ergodic theorem for Markov  (i.e. $L1 =1$) operators
applied to the operator $ \displaystyle \psi (\alpha, z) \to \frac{1}{\varphi
_\alpha (z) } L_{\alpha, f} (\varphi _\alpha \psi (\alpha, . )) (z)$ and the
initial function $ \displaystyle  \frac{1}{\varphi
_\alpha (z) } $.

In order to prove Claim 1, we compute, for a function $\varphi \in L^2$,
$\varphi(z) = \sum _n \gamma _n z^n$, the Fourier coefficients $\gamma '_n$ of
the function $L_{\alpha ,f} \varphi $. We find, for $n\geq 0$:
\begin{eqnarray*}
\gamma ' _n & = &\int \bar z^n L_{\alpha, f} \varphi (z) dz \cr
& = &\int \bar z^n k(\alpha \overline {z}
{\mathstrut f (y)} )
\varphi (y) dydz \cr
&= &\alpha ^n \hat k (-n) \int ({\mathstrut f(y)})^n \varphi (y) dy
\cr
& = &\alpha ^n \hat k(-n) \int \sum _{k\geq 0} {\mathstrut c_k^{(n)}}
y^k \varphi (y) dy \cr
& = &\alpha ^n \hat k(-n)  \sum _{k\geq 0} {\mathstrut c_k^{(n)}}
\gamma _k
\end{eqnarray*}
and, analogously:
$$ \gamma '_{-n} = \alpha ^{-n} \hat k(n)  \sum _{k\geq 0}
{\mathstrut c_k^{(n)}}
\gamma _{-k},$$
where we wrote $(f(z))^n = \sum _{k\geq 0}
c_k^{(n)} z^k, (f(z))^0 = 1$, and $\hat k$ is the Fourier transform of $k$.

Iterating these formulas, we obtain  for a function $\varphi \in L^2$,
$\varphi(z) = \sum _n \gamma _n z^n$, the Fourier coefficients $\gamma
^{(m)}_n$
of the function $L_{\alpha ,f}^m \varphi $:
$$ \gamma _n^{(m)} = \alpha ^n \hat k(-n) \sum _{n_1, \dots , n_m \geq 0} (\Pi
_{s=1}^{m-1} \alpha ^{n_s} \hat k(-n_s)) {\mathstrut c_{n_1}^{(n)}} \;
{\mathstrut c_{n_2}^{(n_1)}} \dots {\mathstrut
c_{n_m}^{(n_{m-1})}} \; \gamma _{n_m} $$ for $n\geq 0$, and
$$ \gamma _{-n}^{(m)} = \alpha ^{-n} \hat k(n) \sum _{n_1, \dots , n_m \geq 0}
(\Pi _{s=1}^{m-1} \alpha ^{-n_s} \hat k(n_s)) {\mathstrut
c_{n_1}^{(n)}} \;
{\mathstrut c_{n_2}^{(n_1)}} \dots {\mathstrut
c_{n_m}^{(n_{m-1})}} \; \gamma _{ -n_m} $$
for $n\leq 0 $.

To get the Fourier coefficients of the bounded continuous function $\int
L_{\alpha, f}^m 1
(z) d\alpha $, we integrate in $\alpha $ the Fourier coefficients of the
bounded continuous functions $ L_{\alpha, f}^m 1 (z) $. In the above sum, all
terms vanish, except the ones with $n+n_1+\dots +n_{m-1} =0 $. Since all $n_i $
have the same sign, the one nonzero integral corresponds to $n= n_1
=\dots = n_{m-1} =0$. Claim 1 follows.
\eproof

\begin{remark}Theorem ~\ref{T=equalmeasures} holds even without
randomization. Suppose that $A$ has determinant equal to $1$ and
let $O$ vary over $SO(2,{\bf R}).$ Then for almost every $O$ the
eigenvalues of $OA$  are either complex with irrational argument
or real and there is one eigenvalue of modulus bigger than one. In
the first case Cesaro sums of the push forward of Lebesgue measure
by $OA_\sharp$ converge to the unique invariant measure of
$OA_\sharp$. In the second case to the Dirac measure supported on
the expanding eigenspace. Call these measures $m_{OA}$, then $\int
m_{OA}dO$ is Lebesgue measure. The proof is the same, but easier.
\end{remark}

\section{Experimental Method} \label{S=Method}
\subsection{Haar Measure on $SO(3)$} \label{SS=Haar}
It is clear that an element of $SO(3)$ is determined
by its axis and angle of rotation.  Here we describe how
to pick axis and angle uniformly with respect to Haar
measure on $SO(3)$.

Let $\Sset^2$ be the usual two sphere with measure $\mu$. Let $\Sset^1$ be
the usual unit circle of angles from $0$ to $2\pi$ given the
probability measure with density function
$(1-\cos\theta)/(2\pi)$. Let $\Sset^2 \times \Sset^1$ be the product
space with the product measure which we denote by $\gamma$. There
is a natural map  $P:\Sset^2 \times \Sset^1 \rightarrow SO(3)$ which maps
a vector $x$ and an angle $\theta$ to the orthogonal
transformation which fixes $x$ and rotates by angle $\theta$
around $x$ according to the right hand rule. The map $P$ sends
$(x,0)$ to the identity in $SO(3)$ for all $x \in \Sset^2$ and is two
to one when $\theta\neq 0$, $(x,\theta)$ and $(-x,-\theta)$ map
to the same point. $P$ maps $\gamma$ to the Haar measure on
$SO(3)$. If we identify $(x,\theta) \sim (-x,-\theta)$ we obtain
$\Sset^2 \times \Sset^1/\sim$ which is a circle bundle over real
projective 2-space. $P$ induces a map $P^\sim:\Sset^2 \times
\Sset^1/\sim\rightarrow SO(3)$.

\begin{prop} $P^\sim$ is one-one off of the zero section
of $\Sset^2 \times \Sset^1/\sim$ and gives a measurable isomorphism
between $(\Sset^2 \times \Sset^1/\sim,\gamma)$ and $(SO(3),\hbox{Haar})$.
\end{prop}

\begin{proof}
That $P^\sim$ is one-one off the zero
section of $\Sset^2 \times \Sset^1/\sim$ is easily verified.

Fix the standard product metric on $\Sset^2\times \Sset^1$,
normalized so that each factor has total volume 1.
Normalized Haar measure on $SO(3)$ is Riemannian volume with
respect to the bi-invariant metric we now describe.
The Lie algebra of $SO(3)$ is the space of anti-symmetric
matrices ${so}(3)$; on this algebra, we put the inner product:
$$<A,B> = \frac{1}{2 c^2} tr(AB^{t}),$$
where $c = 2 \pi^{2/3}$.
An orthonormal basis for ${so}(3)$ is $\{X,Y,Z\}$,
where
$$X =\left(\begin{array}{ccc}
0 & 0 & 0\\
0 & 0 & -c\\
0 & c & 0
\end{array}
\right),\,
Y =\left(\begin{array}{ccc}
0 & 0 & -c\\
0 & 0 & 0\\
c & 0 & 0
\end{array}
\right),$$ and
$$Z =\left(\begin{array}{ccc}
0 & c & 0\\
-c & 0 & 0\\
0 & 0 & 0
\end{array}
\right).$$ (Note $[X,Y] = c Z$, etc). In the bi-invariant metric
induced by this inner product, $SO(3)$ has constant sectional
curvatures, see \cite{DoC}:
$$\lambda = \frac{c^2}{4} = \frac{\|[X,Y]\|^2}{4} =  \frac{\|[X,Z]\|^2}{4}
= \frac{\|[Y,Z]\|^2}{4}.$$  The diameter of $SO(3)$ in this
metric is $r=\pi/c$ and the total volume is
$$\frac{\pi}{\lambda}\left(2r - \frac{\sin( 2\sqrt{\lambda}
r)}{\sqrt{\lambda}}\right) = 1.$$

Let $\rho(x,\theta) d\mu(x) d\theta$ be the pullback
of the volume form on $SO(3)$ to $\Sset^2\times \Sset^1$ under $P$.
To prove that $P^\sim$ is an isomorphism, it suffices
to show that $\rho$ is the density of $\gamma$, that is, to show that:
$$\rho(x,\theta) = (1-\cos\theta)/\pi,$$
for all $(x,\theta)\in \Sset^2\times \Sset^1$.

For any $x,y\in \Sset^2$, if $N=A\times \Theta\subset \Sset^2\times
\Sset^1$
is a product neighborhood of $(x,\theta)$ of Lebesgue measure
$\delta$, then there exists a $g\in SO(3)$ such that
$\hat{N} = gA\times \Theta$ is a neighborhood of $(y,\theta)$ of
Lebesgue measure $\delta$. From the definition of $P$,
it follows that $P(\hat{N}) = gP(N)g^{-1}$.
Since $\rho(x,\theta) d\mu(x) d\theta$
is the pullback of an $SO(3)$-invariant form,
we obtain that for any $x,y\in \Sset^2$, $\theta\in \Sset^1$,
$$\rho(x, \theta) = \rho(y, \theta) =: \rho(\theta).$$

Finally, we compute $\rho(\theta)$.  Since the
geodesics of $SO(3)$ through $I$ are precisely the
one-parameter subgroups,  the image under $P^\sim$ of the curve
$t\mapsto (p, t)$ is a geodesic through $I$ of speed
$1/c$.
It follows that for any $\theta \in (0,\pi)$,
$P^\sim$ sends $\Sset^2\times (0,\theta)/\sim$ diffeomorphically onto
$B_{c^{-1}\theta}  (I) \setminus \{0\}$,
the punctured ball of radius $c^{-1} \theta$ about the identity
in $SO(3)$.
The volume of such a ball is
\begin{eqnarray*}
\int_0^\theta \rho(t)\dd t &=& \hbox{vol}(B_{c^{-1}\theta}(I))\\
&=&\frac{\pi}{\lambda}\left(2c^{-1}\theta - \frac{\sin(
2\sqrt{\lambda}c^{-1}\theta}{\sqrt{\lambda}}\right)\\
& = & \frac{1}{\pi}(\theta - \sin\theta).
\end{eqnarray*}
By the Fundamental Theorem of Calculus,
\begin{equation}
\rho(\theta) = \frac{1}{\pi}(1-\cos(\theta)), \label{e=densrot}
\end{equation}
which completes the proof.
\eproof\end{proof}

We call $P$
or more appropriately $P^\sim$ polar coordinates on $SO(3)$.

\subsection{Computing Random Exponents}

Recall that $R(\eps)=
\int_{P\Sset^2}\log\|Tf_\varepsilon v\|\dd m(v).$ Using the results
of \cite{AvBo} we can reduce the right hand integral to a one
variable integral which we can then evaluate numerically very
accurately. This is how the random Lyapunov exponents are
computed.

\begin{prop} \label{P=formrandom}
$$R(\eps) = \int_0^\frac{1}{2}\log(1 + (2\pi\varepsilon
x(1-x))^2)\dd x .$$
\end{prop}

\noindent{\bf Proof.}
Consider the ``inverse Archimedean projection'' $$\Psi:(\theta,x)\mapsto
(g(x)\cos\theta,
g(x)\sin\theta, x-\frac{1}{2}),$$
where $g(x) = \sqrt{x(1-x)}$.
 This map sends the cylinder
$C=\Sset^1\times[0,1]$ onto the sphere $\Sset^2$ and is area-preserving:
the pullback  $\Psi^\ast d\mu$ is a multiple of
the Lebesgue volume form $d\theta\,dx$ on $C$. The Riemannian
metric on $\Sset^2$ pulls back to the metric:
$$<v,w>_{(\theta,x)} = v^t B(x)^2 w,$$
where
$$B(x)=[D\Psi(x,
\theta)^t D\Psi(x,\theta)]^{1/2} =
\left(\begin{array}{cc}
g(x) & 0\\
0 & (2g(x))^{-1}
\end{array}\right).
$$

Setting $\tilde{f}_\eps = \Psi^{-1}\circ
f_\eps\circ \Psi$, we have that:
$$\tilde{f}_\eps(\theta,x) = (\theta + 2\pi\eps x, x).$$
We next compute:
\begin{eqnarray*}
R(\eps) &=&
\int_{P\Sset^2}\log\|Tf_\varepsilon v\|\dd m(v).\\
& = & \int_{T_1C}\log\|T\tilde{f}_\varepsilon v\|\,\, \Psi_\sharp^\ast \dd
m(v).\\
\end{eqnarray*}
In the second equation, the unit tangent bundle $T_1C$ and the
quantity $\|T\tilde{f}_\varepsilon v\|$ are defined with respect to
the $\Psi$-pullback Riemannian metric on $C$.

If $v\in TC$ is a unit vector with respect to the
pullback metric, then $u = Bv$ is a unit vector with
respect to the Euclidean metric, and
$\|T\tilde{f}_\varepsilon v\| = \|B T\tilde{f}_\eps
B^{-1} u \|_{\hbox{\tiny Eucl.}}$.  Hence we can write:
\begin{eqnarray*}
R(\eps) & = & \frac{1}{2\pi}\int_{(x,\theta)\in C}\int_{u\in \Sset^1} \log
\|B T_{(x,\theta)}\tilde{f}_\eps B^{-1} u\|_{\hbox{\tiny Eucl.}}\, \dd x
\dd\theta \dd u.\\
&=& \frac{1}{2\pi}\int_{(x,\theta)\in C}\int_{u\in \Sset^1} \log
\|\left(\begin{array}{cc}
1 & 4\pi\eps x(1-x)\\
0 & 1
\end{array}\right)
   u\|_{\hbox{\tiny Eucl.}}\dd x \dd\theta \dd u.
\end{eqnarray*}

For $A\in SL(2,{\Rset})$,
\cite{AvBo} show that
$$\int_{u\in \Sset^1} \log \|A u\|_{\hbox{\tiny Eucl.}} \dd u =
\log((\|A\|_{\hbox{\tiny Eucl.}} + \|A\|_{\hbox{\tiny Eucl.}}^{-1})/2).$$
Applying this to $A=\left(\begin{array}{cc}
1 & \alpha \\
0 & 1
\end{array}
\right)$, we obtain that
$$\int_{u\in \Sset^1} \log \|A u\|_{\hbox{\tiny Eucl.}} \dd u =
\frac{1}{2}\log(1 + \alpha^2/4).$$
It follows that:
\begin{eqnarray*}
R(\eps) & = & \int_0^1 \int_{u\in \Sset^1} \log
\|\left(\begin{array}{cc}
1 & 4\pi\eps x(1-x)\\
0 & 1
\end{array}\right)
   u\|_{\hbox{\tiny Eucl.}}\dd u \dd x\\
& = &  \frac{1}{2}\int_0^1  \log(1 + (4\pi\eps x(1-x))^2/4)\dd x\\
& = & \int_0^{1/2}  \log(1 + (2\pi\eps x(1-x))^2)\dd x.
\end{eqnarray*}

\section{Experimental Results} \label{S=Numerics}

In this section we describe several numerical experiments and
their results. First we obtain experimental values of $R(\eps)$
for random maps as introduced in Section~\ref{S=Facts}. We use a
method similar to those we use later to compute estimates for
$\Lambda(\eps)$.  Proposition \ref{P=formrandom} allows to us to
check the accuracy of the computed estimate of $R(\eps)$ against
the precise value of  $R(\eps)$ given there. Then we pass to the
computation of $\Lambda(\eps).$ Three different approaches to the
computation  of $\Lambda(\eps)$ are presented and discussed
which will allow us to obtain accurate enough values to
draw conclusions. Finally, a sample of numerical estimates of
$\Lambda(\eps)$ is shown.

\subsection{The case of random maps} \label{SS=random}

To obtain experimental values of $R(\eps)$ for different $\eps$ we
proceed as in Proposition \ref{P=random} . A random point $x$ in
$\Sset^2$ and a random vector $\xi$ in $T_{1,x}\Sset^2$ is chosen.
A random sequence $g_i$ in  $SO(3)$ is selected and the derivative
of the maps $g_if_\eps$  are applied to the tangent vector $\xi$.
The  rate of increase of the logarithm of $\|T_x f^{(n)}(\xi)\|$
is described, where $f^{(n)}=g_nf_\eps \dots g_1f_\eps.$ The
results are the same with probability 1. For brevity, we refer
loosely to the use of formula (\ref{e=ran1}).

Let us describe the selection of random elements:
\begin{itemize}
\item An initial point can be described in polar coordinates by a
longitude $\lambda_x$ and a latitude $\beta_x$. The value of
$\lambda_x$ is chosen at random in $[0,2\pi]$ with uniform
probability. Concerning $\beta_x$, a random value $z\in[-1,1]$ is
selected with uniform probability and then we let
$\beta_x=\sin^{-1}z.$ This gives the uniform probability for $x\in
\Sset^2.$

\item A tangent vector $\xi\in T_{1,x}\Sset^2$ is generated by
choosing an angle $\psi\in [0,2\pi]$ with uniform probability and
letting  $\xi$ make an angle $\psi$ with the unit tangent vector
to the latitude through $x$ taken in the positive sense. We call
this last vector the horizontal vector at $x$. It is not defined
at the poles, but as the poles have measure zero this is irrelevant
in the current context.

\item A random rotation $g\in SO(3)$ is determined by an axis and
an angle of rotation. As described in subsection \ref{SS=Haar} one
can take these as elements  in $\Sset^2\times \Sset^1.$ It is not
necessary to carry out the identification described there. The
axis is selected just as the point $x\in \Sset^2$ was above. Let
$\theta$ be the rotation angle. To select it according to formula
(\ref{e=densrot}) pick a random value $z\in[0,2\pi]$ with uniform
probability and solve  the equation $z=\theta-\sin(\theta)$ for
$\theta.$ The equation is nothing other than the well known Kepler
equation with eccentricity equal to 1. There are efficient solvers
for it.
\end{itemize}

Then, given initial values of $(x,\xi)$ one can apply formula
(\ref{e=ran1}) to approximate $R(x,\eps)$ by using a finite number
of iterates, $N.$ In turn, to approximate $R(\eps)$ the integral
in formula (\ref{e=ran2}) can be computed using a sample of size
$M$ in $P\Sset^2.$ Let $R_{N,M}(\eps)$ be the value obtained.

This value is compared to the one given by Proposition
\ref{P=formrandom} which has been computed using a Simpson method
with iterative mesh refinement (the values of $R(\eps)$ are shown
in figure XXX). The following limit
approximations are straightforward to derive:
\begin{equation} \label{e=Rlimiteps}
\begin{array}{rclll}
R(\eps)&=&\frac{\pi^2}{15}\eps^2-\frac{2\pi^4}{315}\eps^4+O(\eps^6)
&{\mbox{for}} &\eps\rightarrow 0, \\   \\
R(\eps)&=&\log(2\pi\eps)-2+\frac{1}{2\eps}+O(\frac{\log(\eps)}
{\eps^2}) &{\mbox{for}} &\eps\rightarrow \infty.
\end{array}
\end{equation}
Skipping the $O$ terms one has relative errors less than $0.01$ if
$\eps<0.30$ and $\eps>3.19$, respectively.

Let
\begin{equation}
\Delta_{N,M}(\eps)=R_{N,M}(\eps)-R(\eps). \label{e=diffran}
\end{equation}
We have verified experimentally that $\Delta_{N,M}(\eps)$ has essentially zero
average and a standard deviation like
$$ \sigma_{N,M}(\eps) \approx \frac{\kappa(\eps)}{\sqrt{MN}},$$
provided $N$ is large enough.

Tests have been done for several choices of $N, M$ and $\eps.$
Using $M=10^3$ and $N=10^k,\,k=2,\ldots,6,$ the values of
$\kappa(\eps)$ have been estimated for $\eps$ ranging from
$10^{-1}$ to $10^3.$ There are no significant differences from
$k=4$ on. Figure xxx displays, for different values of
$\eps$, the interval $[-3\kappa(\eps), 3\kappa(\eps)]$ and the
results of single runs (i.e., taking $M=1$) for $N=10^k,\, k=4,
\ldots,9.$ More concretely, the plotted values are the deviations
$\Delta_{N,1}(\eps)$ given by formula (\ref{e=diffran}) multiplied
by $\sqrt{N}.$

These results indicate that $R_{N,M}(\eps)$ and $R(\eps)$ agree to
order $10^{-5}$ taking $NM=10^{10}.$ Further checks have been done
for larger values of $\eps$ (up to $\eps=10^6$) which show no
significant variation of $\kappa(\eps)$ between $\eps=10^3$ and
$\eps=10^6.$

\subsection{Computing the Lyapunov exponent in the deterministic case}
\label{SS=determ}

In principle one can follow a similar scheme to compute
$\Lambda(\eps)$. That is, using (\ref{e=lya1}) with a finite
number of iterates, $N$, an approximation of $\lambda_1(x,g\circ
f_\eps)$ is obtained. Then, (\ref{e=lya2}) is computed using a
Montecarlo method sampling $x\in \Sset^2$ as described with
samples of size $M_p.$ Finally an estimate of $\Lambda(\eps)$ is
computed by applying again a Montecarlo method to (\ref{e=lya3}),
sampling $g\in SO(3)$ as explained above and using samples of size
$M_r.$ In any case the samples are taken according to the
appropriate measures. The total number of iterates of the maps and
their differential is, hence, $M_rM_pN.$ Let us denote by
$\Lambda_{M_r,M_p,N}(\eps)$ a value obtained in this way.

Results of this approach are shown in figure xxx for
different values of $\eps.$ They require some explanation. For a
fixed $g\in SO(3)$ the values of $\lambda_1(x,g\circ f_\eps)$ are
estimated for $M_p$ random values of $x.$ The  standard deviation
of the values of $\lambda_1(x,g\circ f_\eps)$ is then computed.
This value, $\sigma_g,$ depends of the choice of $g.$ Let
$\sigma_{\Sset^2}$ be the average value of $\sigma_g$ when a full
sample of $g\in SO(3)$ is considered. On the other hand all the
$M_rM_p$ determinations of $\lambda_1(x,g\circ f_\eps)$ can be
used to estimate a global standard deviation,
$\sigma_{\mbox{\tiny{total}}}.$

It is clear that $\sigma_{\Sset^2}$ measures the average
dispersion of the maximal Lyapunov exponent when different points
are taken in the phase space. The dispersion depends on the
concrete rotation $g$ taken. Typically, for the $g$ such that
relatively small values of the average $\lambda(g\circ f_\eps)$ of
$\lambda_1(x,g\circ f_\eps)$ on $\Sset^2$ are obtained, it is seen
that $\sigma_g$ is larger. It should also be mentioned that the
errors in the determination of $\lambda_1(x,g\circ f_\eps)$, due
to the finiteness of the number of iterates $N$, also contribute
to this dispersion.

On the other hand there is also a dispersion in $\lambda(g\circ f_\eps)$
when different $g$ are taken. The standard deviation
$\sigma_{\mbox{\tiny{total}}}$ measures the cumulative effect of both
dispersions.

Now let us make several comments on the observed behavior, based on
computations carried out with quite different values of $M_r, M_p$ and
$N.$
\begin{itemize}
\item The estimates of $\Lambda(\eps)$ are close to $R(\eps)$ for
$\eps$ large. For instance, for $\eps>3$ one has that
$|\Lambda_{5000,5000,8192}(\eps)-R(\eps)|$, in the runs done, is
below $0.00325$. The situation is worse for small $\eps$ because
the difference reaches the value $-0.105$ for $\eps=0.55.$ If we
proceed to compute the relative error (r.e.)
$=(\Lambda_{5000,5000,8192}(\eps)- R(\eps))/ R(\eps)$ the
agreement is even worse. For $\eps=1$ one has r.e. $\approx -0.14$
and r.e. $<-0.9$ if $\eps<0.42.$ In these comparisons one should
take into account that the discrepancies also include the errors
done in the estimates of $\Lambda(\eps)$ (see Sections
\ref{SS=improved} and \ref{SS=results}). But the present results
already indicate that the differences for $\eps$ small are not
only due to statistical errors. To make this more evident some
additional computations have used a total number of iterates (for
some selected $\eps$) with $M_rM_pN$ largely exceeding $10^{12}.$

Furthermore it seems also clear that for $\eps<0.3$ there is a
``saturation'' in the behavior of the estimates of
$\Lambda_{M_r,M_p,N}(\eps)$ and of the standard deviations. Indeed, it
can be seen that the trend on the figure changes completely (this is
also the purpose to use logarithmic scales). Systematic errors occur
which completely invalidate the statistical results. To make this more
evident some values, for $\eps$ in the range $[0.2,0.3]$, computed also
in a probabilistic way but with a different estimator (see Section
\ref{SS=improved}) are also shown as dots in the lower left part.
For these computations $M_r=14400,\,M_p=16384$ and $N=16384$ have been
used.

The maximal value $\lambda_{1,{\mbox{\tiny{max}}}}(\eps)$ of
$\lambda_1(x,g\circ f_\eps)$ for $x\in \Sset^2$ and $g\in SO(3)$ is larger
than $R(\eps).$ This requires samples with $M_r, M_p$ large to be
detected if $\eps$ is small. It will be clear from Section
\ref{S=smalleps} and the upper formula in (\ref{e=Rlimiteps}) that
the quotient $q(\eps)=\lambda_{1,{\mbox{\tiny{max}}}}(\eps)/R(\eps)$ is
unbounded if $\eps \to 0.$ If large values of $\eps$ are considered, it
is observed that $q(\eps)$ tends slowly to 1 when $\eps\to\infty.$
In fact it follows from the lower formula in (\ref{e=Rlimiteps}) and
the analysis in Section \ref{S=smalleps} (which is partly valid for any
$\eps$) that the difference $\lambda_{1,{\mbox{\tiny{max}}}}(\eps)-
R(\eps)$ is bounded by $2-\log(2)+O(\eps^{-1}).$ To see differences
close to the bound one has to use very large values of $M_r, M_p$.

\item The estimates of $\sigma_{\Sset^2}$ are mildly sensitive to the
concrete values of $M_r, M_p$ and $N$, provided these values are not too
small, and assuming $0.3<\eps<3.$ For $\eps>30$ a clear dependence with
respect to $N$, of the form $N^{-1/2}$, is seen. From $\eps=3$ to
$\eps=30$ there is a gradual increase in the dependence with respect to
$N.$ For $\eps<0.3$ a tendency towards a behavior of the form $N^{-1}$,
which increases when $\eps$ decreases, is clear.

If extremal values of $\sigma_{\Sset^2}$ are considered when a
sample of $g$ is taken, it is clear that the minimum must be zero.
But there is a significant difference for $\eps<2,$ because the
minimum is already close to zero for samples of moderate size,
while for larger $\eps$ the minimum, which is almost insensitive
to $\eps$ goes to zero slowly when $M_r$ increases. On the other
hand, the maximum of the observed values of $\sigma_{\Sset^2}$
increases until $\eps\approx 10.$ It only stabilizes to a value
with small dependence on $\eps$ for $\eps\approx 100.$

\item The values of $\sigma_{\mbox{\tiny{total}}}$ are much larger than
$\sigma_{\Sset^2}$ for $\eps>1$, while for $\eps<0.3$ they are essentially
equal. In particular they have a small dependence with respect to
$M_r, M_p$ and $N$, if these are not too small and assuming $\eps>0.3.$
\end{itemize}

An analysis of the reasons of the observed behavior is useful because
it helps in three different aspects: a) to understand the different
contributions to the errors in the estimates of $\Lambda(\eps)$; b) to
see the main differences between the cases of random and deterministic
maps; c) to suggest alternative methods to obtain better estimates.
\begin{itemize}
\item In the deterministic case the dynamics on $\Sset^2$ for a
given rotation $g$ is relevant. This dynamics is ``destroyed'' (or
``smoothed'', ``averaged'') in the case of random maps. Hence, the
initial point is irrelevant and the estimates improve in a
probabilistic way, depending on the total number of iterates $MN$,
in the random case. \item Changing $g$ in the deterministic case
produces dramatic changes in the dynamics and, hence, on $\lambda
(g\circ f_\eps).$ This is specially clear for $g$ close to the
identity (axis close to the pole or small rotation) or for
rotations of angle very close to $\pi$ around an axis of small
latitude. The variability of $\lambda (g\circ f_\eps)$ with
respect to $g$ is a major source of dispersion in the results,
specially for large $\eps.$ A standard deviation
$\sigma_{\mbox{\tiny{total}}}$ around $0.4$, mainly due to the
variation of $\lambda(g\circ f_\eps),$ would require sampling with
$M_r$ of the order of $10^8$, at least, to have accurate results.
\item For a fixed $g$, changing $x\in \Sset^2$ has a very mild
effect for large $\eps.$ Despite the possible existence of tiny
islands (see Section \ref{S=smalleps}) the dynamics "looks "
ergodic. A similar behavior has been observed for standard--like
maps in \cite{GLS} and in the case of volume preserving flows it
is seen in \cite{NSS}, where an analysis of the places where the
islands should be expected is carried out before finding them
explicitly.

Figure xxx displays a sample of orbits in $\Sset^2$ for fixed $g$
and
two different values of $\eps$. For $\eps=0.3$ the dynamics is mainly
dominated by an integrable behavior, with many invariant curves and
small chaotic seas, the largest one seen in the front part. This is
persistent with respect to changes in $g.$ The system is even more
integrable (that is, invariant curves fill up a larger area) for most of
the rotations $g\in SO(3).$ For $\eps=2$ (which is not so large!) only
minor islands subsist, and they can even disappear for different $g.$
For values like $\eps=10$ it is hard to see any island unless $g$ is
selected on a set of small measure. In the random case one observes a
uniform distribution of iterates in $\Sset^2$ and the same is essentially
true for large $\eps$ in the deterministic case.

\item On the other hand, for fixed $g$ and small $\eps$ the value
of $\lambda_1(x,g\circ f_\eps)$ depends strongly on $x$. But the
behavior is typically rather sharp. Either one obtains a moderate
value of the order of $\eps^{1.5}$ or it is zero. The smaller the
value of $\eps,$ the larger the measure of the $x$ with exponent
zero, of course. The average value can be very small and despite
the standard deviation $\sigma_{\Sset^2}$ also being small, large
samples with respect to $x$ have to be taken if small relative
errors are desired.

\item The worst point concerning accuracy, especially for small
$\eps$, are the errors in the computation of $\lambda_1(x,g\circ
f_\eps).$ Indeed, for an integrable motion (e.g., $x$ in an
invariant curve) $\lambda_1$ is zero, but the estimates
$\frac{1}{n} \log \|T_xf^n\|$ are, generically, of the order of
$\frac{\log(n)}{n}.$ This implies that the convergence to zero is
slow. For large $\eps$ this effect is relatively not so dramatic,
but oscillations in the behavior of the quotients and different
trends can be expected.

\end{itemize}

\subsection{Improved procedures} \label{SS=improved}

Several alternative procedures have been used in previous computations
to determine the maximal Lyapunov exponent of a given map, averaged
on the phase space:
\begin{itemize}
\item [1)] If the dynamics has a uniform hyperbolicity but with
superimposed strong periodic or quasi\-periodic oscillations, the
following strategy has been used in \cite{BS}. It consists in detecting,
by an iterative procedure, an upper envelope of the plot of the
quotients $\frac{1}{n} \log \|T_xf^n\|$ as a function of $n.$ Then, and
after skipping a transient regime, one fits a function of the form
$\alpha+\beta/n$ to the envelope. The value of $\alpha$ is a good
estimator for $\lambda_1(x,f)$ and, as the system in \cite{BS} is a skew
product with linear action on the fibers, the value of $x$ is
irrelevant.
\item [2)] If the values of $\lambda_1(x,f)$ depend strongly on $x$, it
is possible to divide the phase space in pixels of a given size (in
general, $d$-dimensional pixels) and start the computations at a point
in each pixel. However, if the number of pixels is large and the system
depends on several additional parameters, the method can be prohibitive
from a computational point of view. Then, together with each initial
pixel one considers all the pixels visited by the orbit. The estimated
Lyapunov exponent is assigned to all of them. One requires each pixel
to be visited a minimal number of times (in case of need one takes
several initial points in the pixel) and an averaged Lyapunov exponent
is assigned to the pixel. Later on this is averaged over the full
phase space.

This method has been used in \cite{SS} to study the classical
Hill's problem and how the degree of chaos behaves with respect to
the energy.

One should also take into account the stickiness of invariant curves.
An initial point in a chaotic sea can remain close to an island for a
large number of iterates. Hence, it is a good strategy to take
a larger number of initial points even if one has to decrease the
number of iterates for each one, provided this number is not too small.
Furthermore, the local slope of $\log \|T_xf^n\|$ can have quite
different trends if the number of iterates is large. Statistically
this is not a problem because the interesting magnitude is the average
behavior.
\item [3)] The previous method still suffers from slow convergence of
the quotients to $\lambda_1(x,f)$. An alternative method has been used
in the context of flows with applications to galactic potentials in
\cite{CS} and later on extended to discrete transformations in
\cite{CGS}, where references to other applications can be found. It is
mainly intended to discriminate between regular and chaotic motion
(that is, to decide if one can accept $\lambda_1(x,f)=0$ or not), but it
 also supplies an estimate of the Lyapunov exponent.

Given a map $f$, an initial point $x$ on a manifold $\CM$ and a
random vector $\xi\in T_{1,x}\CM$, let $\xi_0=\xi$ and define
$\xi_k=(T_{f^{k-1}(x)}f)(\xi_{k-1}).$ For fixed integers $m$ and $n$ and
$N>0$ let
$$Y_{m,n}(N)=N^n\sum_{k=1}^N \log\left(\frac{\|\xi_k\|}{\|\xi_{k-1}\|}
\right) k^m,\qquad\bar{Y}_{m,n}(N)=\sum_{k=1}^N Y_{m,n}(k). $$
Then, for a chaotic orbit the estimator
\begin{equation}
\hat{Y}_{m,n}(N)=(m+1)(m+n+2)\frac{\bar{Y}_{m,n}(N)} {N^{n+m+2}}
\label{e=megno}
\end{equation}
tends to $\lambda_1(x,f)$, while for a generic regular orbit it
behaves like $\frac{(m+1)(m+n+2)}{m(m+n+1)}\frac{1}{N}.$ The basic
idea is to average the exponential rate of increase of the length
of $\xi_k$ so that the transience has small relevance and to smooth
out the irregularities of the quotients. Hence, it is a measure of
the mean exponential growth of nearby orbits (MEGNO) and depends
on the couple $(m,n).$ Suitable values (according to numerical
experience) are $m=2,\,n=0,$ and then it is denoted as MEGNO20.

For $(m,n)=(2,0)$ a slightly better estimator for $\lambda_1(x,f)$ is
obtained by using $12\frac{\bar{Y}_{2,0}(N)} {N^4+4N^3+5N^2}.$
Furthermore, when this method is used with these $(m,n)$, one can check
for a behavior of the form $\frac{2}{N}$ to decide $\lambda_1(x,f)
\approx 0.$ Typically $\hat{Y}_{2,0}(N)-\frac{2}{N}=O(N^{-2})$ for
regular orbits.

Given a maximal number of iterates $N_{\mbox{\tiny{max}}}$ to be
used in the estimates, an additional question is whether it can be
better to use another value $N<N_{\mbox{\tiny{max}}}$ as a better
choice to estimate $\lambda_1(x,f)$. In \cite{CGS} a ``right
stop'' criterion is introduced. It is specially relevant if the
orbit is close to be regular, to prevent an overestimate of
$\lambda_1(x,f)$, but it has not been used in the present
computations.

Figure xxx illustrates the different behavior of MEGNO20
and the quotients in formula (\ref{e=lya1}) in several cases. Details
on the parameters used for the plots are as follows.

\end{itemize}

On the upper row, left plot, for $\eps=0.3$ a rotation of angle
$\theta$ and axis of latitude $\beta$ with $\frac{\beta}{2\pi}=\frac{3}
{25}g,\,\frac{\theta}{2\pi}=\frac{18}{25}g$ where $g=(\sqrt{5}-1)/2$ has
been selected. Random initial conditions are chosen. After a transient
of 512 iterates, estimates of $\lambda_1(x,f)$ are produced and plotted
for the next 2048 iterates. Solid lines (the lower ones) correspond to
the estimates using MEGNO20, while discontinuous lines are produced by
the classical formula. The middle part displays a similar plot for
different initial conditions and the right part is similar to the middle
one but for $\eps=10.$

On the lower row estimates of the largest Lyapunov exponent for
random initial conditions are shown. On the left part $\eps=0.3$
and $\theta$ and $\beta$ as before have been used, and 512 random
initial points are plotted. Solid (resp. discontinuous) lines
correspond to MEGNO20 (resp. classical) estimators. Transient and
number of iterates are as before. The left (middle) plot in the
upper row corresponds to number 4 (43) of these points. The middle
part is similar but for $\eps=3.$ For completeness estimates in
the random case using finite $n$ in (\ref{e=lya1}) are also shown
on the right for $\eps=0.3\,.$ Each iterate uses a random element
in $SO(3)$ with density $\hat{\rho}(\beta,\theta)=
\cos(\beta)(1-\cos(\theta))/(2\pi).$ The solid line shows the
value of $R(0.3)\approx 0.0547518\,.$

See also \cite{S-tak} for additional methods and applications.

Due to the good properties of the procedure,  estimates of
$\lambda_1(x,g\circ f_\eps)$ have been computed using procedure 3)
above. For the integrations in $\Sset^2$ a Montecarlo method has
been used. This is good enough for large $\eps$, because of the
mild dependence of $\lambda_1(x,g\circ f_\eps)$ with respect to
$x$ (for, say, $\eps>1$) for most of the $x\in \Sset^2$ and most
of the $g\in SO(3).$ For $\eps$ small and especially if
$\eps<0.3,$ a method such as the one presented in 2) would be
suitable, but there are additional problems, due to the smallness
of the Lyapunov exponents, to be discussed later. Furthermore it
will turn out that it is relevant to compute $\Lambda(\eps)$ with
a small relative error for large $\eps,$ to allow for a careful
comparison with $R(\eps).$ But for $\eps$ small it will be clear
from the results, even those obtained with a moderate accuracy,
that $\Lambda(\eps)$ is far below $R(\eps).$ In any case, it seems
that numerical estimates of $\Lambda(\eps)$ for $\eps<0.3$ with
small relative error require an enormous computational effort.

Finally, for the integrations in $SO(3)$ and taking into account
the large standard deviation found for moderate and large values
of $\eps$, it has been found more convenient to use numerical
quadrature formulas based on a grid of points. More concretely, a
product Simpson method has been used with respect to the latitude
$\beta_g$ of the axis of rotation and the rotated angle
$\theta_g.$ The longitude of the axis is irrelevant: changes in
this longitude are equivalent to changes in the longitude
$\lambda_x$ of $x\in \Sset^2.$ Using a grid with $\theta_g\in
[0,2\pi],\,\beta_g\in[0,\pi/2]$ requires that the estimate of
$\lambda(g\circ f_\eps)$ be multiplied by the factor
$\cos(\beta_g)(1-\cos(\theta_g)).$

Initial estimates for the results shown in the next Section use,
for the elements in $SO(3)$ the following data: $\lambda=0$,
$(\theta,\beta)$ on a grid of $N_g\times N_g$ with $N_g=200$. Then
a sample of $N_p=1600$ random initial points and tangent vectors
in $P\Sset^1$ and $M=8000$ iterates are used. The programs have
been done in such a way that beyond the estimates for these
values, also estimates using grids with $N_g=100,\,50$ or using
samples with $N_p=800,\,400,\,200,\,100$ and doing a number of
iterates equal to $M=4000,\,2000,\,1000,\,500$ are computed. This
allows for a check on the internal consistency of the results.

It turns out that the use of different grids in $SO(3)$ stabilizes
quickly. Concerning the dependence with respect to $N_p$ and $M$,
it is clearly seen that there is no need for very large values of
$N_p$ except in the case of small $\eps$ and one is interested in
having small relative error. The dependence with respect to $M$ is
clearly of the form ctant/$M.$ Hence, extrapolations with respect
to $M$ have been used. The initial estimates allow for a fine
tuning of the most suitable values for the grid, $N_p$ and $M$.
For instance, assuming that one can accept a total of $2^{41}$
iterates (for every value of $\eps$), for large $\eps$ a typical
choice is $N_g=2^8,\,N_p=2^9,\,M=2^{16},$ while for small $\eps$
it is $N_g=2^7,\,N_p=2^{10},\,M=2^{17}.$ Even with this large $M$
the results start to be not very good if $\eps<0.2\,.$

Figure xxx shows 3D views of $h(\theta_g,\beta_g,\eps)=
\lambda(g\circ f_\eps) \cos(\beta_g)(1-\cos(\theta_g))\frac{\pi}{2}$ as
a function of $(\theta_g,\beta_g)$ for different values of $\eps.$ Level
lines of these surfaces are displayed in figure xxx. The plots
give a good evidence of the smooth behavior of $h(\theta_g,\beta_g,
\eps)$ for moderate and large values of $\eps$, and how the behavior
becomes wilder, with sharp changes for small $\eps.$ It is clear that
the results are the same if $(\theta,\beta)$ is replaced by
$(2\pi-\theta,-\beta)$. Furthermore, note that for integer $\eps$ the
symmetries of $g\circ f_\eps$ imply that the results should be the same
if $\theta$ is changed to $2\pi-\theta$, as clearly seen in the first
three plots. For $\eps=0.5$ the lack of symmetry $\theta \leftrightarrow
2\pi-\theta$ is clear, but for large non--integer $\eps$ this lack of
symmetry is harder to detect.

\subsection{A sample of results}\label{SS=results}

The results of applying the methodology just described are shown in
figure xxx. Typical values for the number of iterates, initial
data and grid have been given before.

On the top plots general views of $R(\eps)$ and estimates of $\Lambda
(\eps)$, to be denoted as $\Lambda(\eps)_{\mbox{\tiny{num}}}$, can be
seen. In particular, for large and small $\eps$ it is easy to check the
limit behavior of $R(\eps)$ predicted by (\ref{e=Rlimiteps}). On top
left only values $\eps\le 10$ are shown. On this scale no differences
can be seen between $R(\eps)$ and $\Lambda(\eps)_{\mbox{\tiny{num}}}$
for $\eps>3.$ For small $\eps$ the differences are clear and they are
quite dramatic for $\eps<0.4.$

On the bottom left plot the tiny differences
$\Lambda(\eps)_{\mbox{\tiny{num}}}-R(\eps)$ are displayed. It
seems that they tend to 0 as $\eps\to\infty,$ in agreement with
the second possibility in Section \ref{S=Heuristics}. Finally the
behavior of $\Lambda(\eps)_{\mbox{\tiny{num}}}$ for $\eps$ small
is seen in detail on the bottom right part. A logarithmic scale
has been used to reveal that
$\log(\Lambda(\eps)_{\mbox{\tiny{num}}})$ is dominated by a
function of the form $-c/\eps$ for some $c>0.$ Due to the
smallness of the estimates and to the fact that it would require
an enormous effort to estimate $\Lambda(\eps)$ for $\eps$ close to
$0.1$ (unless other methods, from deterministic analysis, are
used), it is not completely clear what the correct behavior is.
Even with a reduced set of data a
fit of the values obtained for $0.16\le\eps\le 0.3$ gives a result of the
form
\begin{equation}
\log(\Lambda(\eps)_{\mbox{\tiny{num}}})\approx 2.45-\frac{3.16}{\eps}
\label{e=estimexp}
\end{equation}
which must be taken with caution, but seems to give the correct
trend. This suggests that the inequality in Question~\ref{Q6}
is not satisfied for small $\eps$, which will be confirmed
theoretically in
Section \ref{S=smalleps}. This fact is not a surprise, because
similar facts occur in generic analytic families of
area-preserving diffeomorphisms. The smallness of the Lyapunov
exponent is related to the area of the chaotic seas which in turn
is related to the splitting of the separatrices of fixed and
periodic point for maps close to the identity. See \cite{FS},
where general upper estimates can be found. In fact, as with many
a priori exponentially small upper bounds, this result can also be
obtained as a corollary of averaging theory for analytic systems,
see \cite{N} and also Section \ref{SS=pertwist}.

Finally it should be mentioned that some computations have been done
for large values of $\eps$ (up to $10^6$). Due to the strong chaotic
properties it is enough to take small values (say $N_p=256$) of the
number of initial points in $\Sset^2$. But the grid in the parameters
$(\theta_g,\beta_g)$ has to contain more points. Typical values of $N_g$
to have a good determination of $\Lambda(\eps)_{\mbox{\tiny{num}}}$ are
$2^9$ and $2^{10}$. The results confirm what is seen in the left lower part
of figure xxx, that is, $\Lambda(\eps)_{\mbox{\tiny{num}}}>R(\eps)$
and the difference goes to zero slowly.
\section{The case of small $\eps$} \label{S=smalleps}
As it is clear that the greatest problems occur for small $\eps$,
it is worth it to carry out a preliminary analysis of the dynamics
in this case. The first item to be studied is the location and
stability of fixed points. This can be carried out, with the same
effort, for any $\eps$. Furthermore this allows us to see how
bifurcations give rise to new elliptic fixed points with the
corresponding creation of islands. Later on the global behavior of
$g\circ f_\eps$ on $\Sset^2$ is discussed. In what follows it is
assumed that $\eps>0.$

\subsection{Fixed points and their stability}\label{SS=fixed}
To look for fixed points of $g\circ f_\eps$ it is enough to consider
axes of rotation with zero longitude. Let $(\beta,\theta)$ be the
latitude of the axis and the angle of the rotation. As follows from
Section \ref{SS=Haar}, it is not restrictive to assume $\beta\in
[0,\pi/2].$ Then a fixed point $A$ is mapped by $f_\eps$ to a point $A'$
which by $g$ returns to $A$. Let $b$ be the latitude of $A$. It is clear
that $A$ and $A'$ must have symmetric longitudes, $-\delta(b)$ and
$\delta(b)$, respectively, where $\delta(b)=\frac{\pi}{2}\eps
(1+\sin(b)).$ It is easy to derive the condition for the fixed points
\begin{equation} \label{e=fixed}
\sin\left(\frac{\theta}{2}\right)\sin(\beta)\cos(b)\cos(\delta(b))\!-\!
\sin\left(\frac{\theta}{2}\right)\cos(\beta)\sin(b)\!+\!
\cos\left(\frac{\theta}{2}\right)\cos(b)\sin(\delta(b))\!=\!0.
\end{equation}
For $\eps$ small one has $\sin(\delta(b))=O(\eps)$ and $\cos(\delta(b))
=1-O(\eps^2).$ For the analysis of this case it is more convenient to
write (\ref{e=fixed}) in the form
\begin{eqnarray} \label{e=fixedeps}
\sin\left(\frac{\theta}{2}\right)\sin(\beta-b)&-&
\eps^2\frac{\pi^2}{4}
\sin\left(\frac{\theta}{2}\right)\sin(\beta)\cos(b)
(1+\sin(b))^2 \nonumber \\
&+& \eps\frac{\pi}{2}\cos\left(\frac{\theta}{2}\right)\cos(b)
(1\!+\!\sin(b))+ O(\eps^3)=0.
\end{eqnarray}

 From (\ref{e=fixedeps}) it is clear that if $\sin\left(\frac{\theta}{2}
\right)=O(1),$ that is, $\theta$ is not too close to 0 from the positive
or the negative side, then one has solutions for $b$ of the form
$$ b=\beta+\gamma\eps+O(\eps^2) \qquad {\mbox{or}} \qquad
b=\pi+\beta+\gamma\eps+O(\eps^2),$$
where $\gamma$ is independent of $\eps.$ The value for $\gamma$ is
given by
$$ \gamma=\frac{\pi}{2} \cos(\beta)(1\pm\sin(\beta))/
\tan\left(\frac{\theta}{2}\right) ,$$
where the $+$ sign is used in the first case and the $-$ sign in the
second. In both cases $b$ is close to either $\beta$ or $\beta+\pi$ and
an important thing is that there are exactly two fixed points for
$g\circ f_\eps.$

Otherwise one can write $\theta=m\pi\eps,\,m>0$ or $\theta=2\pi+
m\pi\eps,\,m<0.$ The dominant terms in (\ref{e=fixedeps}) become
in both cases
\begin{equation} \label{e=limiteq}
m\sin(b-\beta)-\cos(b)(1+\sin(b))=0.
\end{equation}
Equation (\ref{e=limiteq}) has to be seen as an equation for $b$
depending on $\beta$ and $m$ (which accounts for $\theta$). As it
has zero average it should have at least two different zeros. To look
for more solutions it is relevant to compute the lines (in $(\beta,m)$)
where double zeros occur. The angle $b$ can be used to parameterize
these lines. They are represented as
\begin{equation}\label{e=doublezer}
\begin{array}{rcl}
m^2 & = & 2+\sin^3(b)-\frac{3}{4}\sin^2(2b),\\ \\
\beta&=&b-\arg\left[\sin(b)-\cos(2b)-\sqrt{-1}\left(\cos(b)+
\frac{1}{2}\sin(2b)\right)\right].
\end{array}
\end{equation}
As $\beta\in[0,\pi/2],$ inspection of (\ref{e=limiteq}) shows that no
double zeros can occur in the case $m>0.$ Hence the value of $m$ is
confined to $[-2,0].$ The bounds on $\beta$ also imply that the
parameter $b$ in (\ref{e=doublezer}) has to be in $[\pi/2,3\pi/2].$ It
is elementary to discuss the behavior of $(m,\beta)$ as a function of
$b$. It is better seen by looking at figure xxx. It can just be
said that two curves of double zeros appear with $\beta\in[\pi/4,
\pi/2].$ They meet at $m=-\sqrt{2},\,\beta=\pi/4,$ where a triple zero
appears. Between both curves there are exactly four zeros.

When additional powers of $\eps$ are included, a routine
application of the Implicit Function Theorem permits us to
conclude the same behavior for the full equation (\ref{e=fixed}).
It should be noted that $\beta=\pi/2$ has to be excluded from the
previous analysis: in that case the axis of rotation is also the
axis of the twist.

The preceding analysis can be summarized as
\begin{prop}\label{P=zeroseps}
For $\eps$ small enough and any $g\in SO(3)$ there are always at least
two fixed points of $g\circ f_\eps.$ Bifurcations to exactly four fixed
point appear for any longitude of the rotation axis and for latitude of
the axis and angle of rotation $(\beta,\theta=\pi(2-m\eps+O(\eps^2)))$
along two lines described by formulae {\rm (\ref{e=doublezer})} when the
parameter $b$ ranges in $(\pi/2,\pi)$ and $(\pi,3\pi/2)$, respectively.
\end{prop}

To discuss bifurcations of the fixed points for general values of $\eps$
is an elementary but cumbersome task. As an illustration the case of
bifurcations appearing on $\beta=0$ is presented. Then (\ref{e=fixed})
reduces to
\begin{equation}\label{e=fixedbet0}
\cos\left(\frac{\theta}{2}\right)\cos(b)\sin(\delta(b))-
\sin\left(\frac{\theta}{2}\right)\sin(b)=0
\end{equation}
and the condition for a double root becomes
\begin{equation}\label{e=doublegen}
\frac{2}{\pi\eps}\tan\left(\frac{\pi\eps}{2}(1+\sin(b)\right)=\sin(b)
\cos^2(b).
\end{equation}
The degenerate cases $|b|=\pi/2$ must be excluded in
(\ref{e=doublegen}). It is immediate that new double fixed points
appear on $\Sset^2$ with $\beta=0$ if and only if $\eps$ is a
positive integer. The number of double fixed points with $\beta=0$
(and some $\theta$) increases with $\eps.$ Also from
(\ref{e=fixedbet0}) it follows that new zeros appear near
$\theta=0,$ one for $\theta>0$ and the other for $\theta<0.$ These
zeros move towards $\theta=\pi$ without ever reaching it. So, it
is a simple matter to state how many fixed points exist for
$\beta=0$ (except at the bifurcation values of $\eps$): there are
at most $2(1+E(\eps))$, where $E$ denotes the integer part of
$\eps.$ For a given non-integer $\eps$ there are always values of
$\theta$ such that this number $2(1+E(\eps))$ is the exact number
of fixed points. This has an elementary dynamical interpretation:
new fixed points emanate from the north pole of $\Sset^2$ when the
rotation number of $f_\eps$ at the north pole (defined by
continuity) passes through 0 (mod 1).

To study the stability of the fixed points we recall that they are
generically elliptic (eigenvalues $\mu$ in $\Sset^1\setminus\{\pm
1\}$), hyperbolic (real positive eigenvalues) and hyperbolic with
reflection (real negative eigenvalues). Let $E,H,R$ denote the
number of fixed points of each type. Euler--Poincar\'e formula
gives $E-H+R=2$ (for simple fixed points). At the creation of new
fixed points $E$ and $H$ increase by 1. When double eigenvalues
are equal to $-1$ then, generically, $E$ decreases by 1 and $R$
increases by 1.

An analytic discussion on the stability of the fixed points is
elementary (at least for small $\eps$) but cumbersome. It is worth
mentioning  that, for any $\eps$ the maximum eigenvalue at a fixed
point is achieved on $b=0$ and has the expression
\begin{equation}\label{e=maxeigen}
\mu_{\mbox{\tiny{max}}}=\frac{\pi\eps}{2}+\left(1+\left(\frac{\pi\eps}
{2}\right)^2\right)^{1/2}.
\end{equation}

A sample of illustrations is shown in figure xxx having
$\theta/2\pi$ as horizontal variable and $\beta/\pi$ as vertical one.
A region containing $i$ (resp. $j,k$) fixed points of type $E$ (resp.
$H,R$) is denoted as $R^kH^jE^i.$ On the top left plot and for
$\eps=0.1$ the two curves on the upper part of the plot
are the curves of double zeros given by (\ref{e=doublezer}). Only in
the region bounded by them there are 4 fixed points; the code is
$H^1E^3.$ The codes for the black, dark grey and light grey regions
are, respectively, $R^2,\,E^2$ and $R^1E^1.$ On the top right plot,
the value of $\eps$ is $\approx 3.456789.$ The
region containing the point $(0.5,0)$ has exactly 2 fixed points while
the regions which contact with this one through arcs have 4. The darker
region has 6 and the small region near the upper right corner has 10.
The regions around $\theta=0$ have 8 fixed points. The solid lines
give the location of all bifurcations and changes of stability.

In the bottom left plot, computed for $\eps\approx 9.876543$, all
the lines of bifurcation or change of stability are plotted. The number
of fixed points, NFP, in the major regions is shown. The typical
transitions are as follow: Consider, for instance, a passage from
NFP$=16$ to NFP$=18$ near $\beta=0$ with increasing $\theta.$ First a
line of creation of an elliptic and hyperbolic points is found. One
passes from a code $R^9H^7$ to $R^9H^8E^1.$ This is followed by a change
of stability by passing to $R^{10}H^8.$ Later on, inside the region with
NFP$=20$, the points of $R$ type become again of type $E$. So the code
passes, in the different changes, from $R^{11}H^9$ to $H^9E^{11}.$

Finally, in the bottom right plot, for $\eps=0.5$, regions similar to
the case $\eps=0.1$ can be seen, with a different configuration. The
level lines of figure xxx are also displayed. It is checked that
the highest levels correspond to domains where the map has exactly one
elliptic and one hyperbolic fixed points. This fact is also present
for smaller values of $\eps.$

\subsection{The maps as perturbed twists}\label{SS=pertwist}
The object of interest is the global dynamics of $g\circ f_\eps$
in $\Sset^2$. To this end it is convenient to write these maps in
a slightly different, but equivalent, way. For this study
$\Sset^2$ will be taken as the sphere of radius 1 centered at the
origin. Instead of considering $f_\eps$ as a twist around the
$z$-axis of angle $\pi\eps(1+z)$, it will be taken as a twist
around the axis of zero longitude and latitude $\beta$ with angle
of rotation around this axis of a point of coordinates $(x,y,z)$
equal to $\pi\eps(1+x\cos(\beta)+z\sin(\beta)).$ Then, the
rotation $g$ is simply a rotation of angle $\theta$ around the
$z$-axis, to be denoted by $R^{(z)}_\theta.$ Up to the
substitution of $\pi/2-\beta$ for $\beta$, the relative positions
of the axes in this formulation and in the previous one are
equivalent.

It is instructive to first consider the case
$\theta=2\pi\frac{p}{q},$ where $p,q$ are coprime integers. Let us
introduce $\delta=\pi(1+ x\cos(\beta)+z\sin(\beta)).$ Then
\begin{equation}\label{e=twistonly}
f_\eps\left(\begin{array}{c}x\\y\\z\end{array}\right) =
\left(\begin{array}{c}x\\y\\z\end{array}\right)+\eps\delta
\left(\begin{array}{c}-y\sin(\beta)\\x\sin(\beta)-z\cos(\beta)\\
y\cos(\beta)\end{array}\right) + O(\eps^2).
\end{equation}
The next step is the computation of the map
$M_{q,\theta,\beta,\eps}:= (R^{(z)}_\theta\circ f_\eps)^q,$ the
parameter $\beta$ being the latitude of the axis of $f_\eps.$ It
is clear that at order zero in $\eps$ one has
$M_{q,\theta,\beta,0}=\mbox{Id}.$ An elementary computation using
formula (\ref{e=twistonly}) for $f_\eps$ and the expression of
$\delta$ as a function of $x,z$ and $\beta,$ gives
\begin{equation}\label{e=themap}
M_{q,\theta,\beta,\eps}\left(\begin{array}{c}x\\y\\z\end{array}\right) =
\left(\begin{array}{c}x\\y\\z\end{array}\right)+
R^{(z)}_\gamma\left(\begin{array}{c}x\\y\\z\end{array}\right)+O(\eps^2),
\end{equation}
where $R^{(z)}_\gamma$ is now a rotation around the $z$-axis in each one
of the horizontal planes with angle of rotation depending on $z$ as
follows
\begin{equation}\label{e=gamma}
\gamma=\pi\eps q\left(\sin(\beta)+zP_2(\sin(\beta))\right),
\end{equation}
where $P_2$ denotes the second Legendre polynomial
($P_2(w)=\frac{3}{2}w^2- \frac{1}{2}$).

This result tells us that the rotation $R^{(z)}_\theta$ averages
 the effect of the map $f_\eps$ in a good way. Let us remark that
in (\ref{e=themap}) it is assumed $p/q$ fixed and $\eps$
sufficiently small. From (\ref{e=gamma}) it follows that the angle
$\gamma$ is still small provided that $\eps q$ is small. In the
trivial case $p=0,q=1$ the rotation is the identity and then the
twist can be be taken with $\beta=\pi/2,$ recovering in formula
(\ref{e=gamma}) the angle rotated in the twist.

To pass to the general case for $\theta$ one needs a preliminary lemma.

\begin{lemma}\label{l=intervals}
Let $\rho\in (0,1)$ and $N\in\Nset.$ Let
$$S_{\rho,N}=\bigcup_{1\le q\le N,\,0\le p\le
q,\,(p,q)=1}\left[\frac{p-\rho}{q},
\frac{p+\rho}{q}\right].$$
Then, if $N+1\ge \rho^{-1}$ one has $[0,1]\subset S_{\rho,N}.$
\end{lemma}

\begin{proof} Let $\alpha\in[0,1].$ If $\alpha=r/s\in\Qset$ with $(r,s)=1$
and
$s\le N,$ there is nothing to prove. Otherwise consider the approximants to
$\alpha$ given by the continued fraction algorithm.  Assume that
$\frac{p_1}{q_1}$ and $\frac{p_2}{q_2}$ are consecutive approximants with
$q_1\le N$ and $q_2 > N.$ Then
$$ \big|\alpha-\frac{p_1}{q_1}\big|\le \frac{1}{q_1q_2}\le
\frac{1}{q_1(N+1)} \le\frac{\rho}{q_1}.  $$
\!\!\eproof
\end{proof}

\bigskip

Hence $\frac{\theta}{2\pi}$ can be written as
$\frac{p}{q}+\frac{\mu}{q}$, where
$|\mu|\le\rho,\,q+1\le\rho^{-1},$ where $\rho$ is not specified
for the moment. Hence one can represent the map $g\circ f_\eps$ as
something similar to the previous case, that is, a rotation whose
angle is a rational multiple of $2\pi$, composed with a map close
to the identity, by writing
$$ g\circ f_\eps= R^{(z)}_{2\pi p/q}\circ R^{(z)}_{2\pi\mu/q}f_\eps.$$
By a direct computation one obtains  expression
(\ref{e=themap}) again with the following modifications:
\begin{itemize}
\item If the rational which approximates $\frac{\theta}{2\pi}$ is 0, then
$q=1$
and there is no average, so that we keep the map $g\circ f_\eps,$
\item The rotation is now $\gamma=\pi\eps q\left(\sin(\beta)+
zP_2(\sin(\beta))\right)+2\pi\mu$,
\item The error terms are, uniformly in $\theta$, of the form
$O(\eps^2 q+\frac{\rho^2}{q}).$
\end{itemize}
If one takes $\rho=\eps^{2/3}$ then the maps $(g\circ f_\eps)^q$
are, in all cases, $\eps^{1/3}$--close to the identity and the
error terms are at most $\eps^{4/3}$. Note that besides the choice
$\rho=\eps^{2/3}$ there are other possibilities, but
$\rho=\eps^{2/3}$ is good enough to prove Corollary \ref{C=corol}.
Finally, it is clear that (\ref{e=gamma}), or the modification
just mentioned adding $2\pi\mu$, is a twist except for $\beta=
\beta_{\mbox{\tiny{crit}}}$ such that
$P_2(\sin(\beta_{\mbox{\tiny{crit}}}))=0$
($\beta_{\mbox{\tiny{crit}}}=\sin^{-1}(1/\sqrt{3})$).

To summarize we state the following
\begin{prop}\label{P=prefinal}
If $\theta$ is not $\eps^{2/3}$--close to zero and $\beta\neq
\beta_{\mbox{\tiny{crit}}}$ the maps $g\circ f_\eps$ for $\eps$ small
enough, have a power which is $\eps^{1/3}$--close to the identity. This
power satisfies a twist condition of order at least $\eps
P_2(\sin(\beta))$.
\end{prop}

\begin{remark}\label{Rem1}
Rotations $g$ with small $\theta$ are irrelevant, for the present purpose,
due
to the Haar measure in $SO(3)$, and a single exceptional case (non-twist)
is
also unimportant. This will be seen later in detail.
\end{remark}

\begin{remark}\label{Rem2}
In the case of small $\theta$ it is still possible to show that
$g\circ f_\eps$ produces a twist effect on each meridian in
$\Sset^2.$ A problem which appears though is that the angle
rotated by the different points can pass through an extremum,
losing in this way the twist property. In fact this is not so
important because the existence of invariant curves when the twist
condition is lost at some point has been established in
{\rm \cite{S-nt}}. But this refinement is not necessary in the present
context.
\end{remark}

\begin{theorem}\label{T=final}
With the possible exclusion of an open set ${\cal B}$ in $SO(3)$ of
small measure, there exist $\eps_0$ such that for $\eps<\eps_0$ the maps
$g\circ f_\eps$ have a dynamics exponentially close to an integrable
flow in $\Sset^2.$
\end{theorem}
\begin{proof}
The proof is divided into  steps.
\begin{enumerate}
\item The maps $M_{q,\theta,\beta,\eps}$, being a power of $g\circ
f_\eps$,
have the same dynamics as $g\circ f_\eps$. In all cases (including $q=1$
and the
exceptional value of $\beta$) they are $\eps^{1/3}$--close to the identity.
Hence there exists a suspension given by the flow of a $1$--periodic
vector
field in $\Sset^2$ such that the time--$1$ map associated to this flow
coincides with $M_{q,\theta,\beta,\eps}$. The vector field is ``slow'' (of
the order of $\eps^{1/3}$) and the dominant terms do not depend on time.
See \cite{BRS} for details and an explicit construction.

It is relevant to note that the vector field is analytic with
respect to the phase space variables (that is, the points in
$\Sset^2$) while the dependence in $\eps$ is discontinuous in
$SO(3)$ (moving $\theta\in[0,2\pi]$ changes the value of $q$), but
the relevant thing is that it is {\em bounded} in $\eps$.
Furthermore the dependence with respect to $t$ can be made of
class ${\mathcal C}^r$ for any $r>0,$ but continuity in $t$ is
sufficient for what follows. Furthermore the vector field is
Hamiltonian. \item The next step is to ``average'' the vector field
with respect to $t$.  This is the content of Neishtadt's theorem
\cite{N}. See \cite{S-fast} for a detailed proof. As a consequence
the vector field can be written as an autonomous part and a
remainder which is exponentially small in the current small
parameter; that is, the remainder is bounded by
$\exp(-c\eps^{-1/3})$ for some $c>0$. Furthermore, the averaged
vector field is still Hamiltonian (see \cite{SV} for a sketch of
the proof). \item As the averaged system is a Hamiltonian in
$\Sset^2$, it is integrable and, hence, foliated by invariant
curves except on the separatrices, which are a set of zero
measure. Most of the invariant curves subsist as a consequence of
Moser's twist theorem. To this end one should have that the
perturbation is small compared with the twist condition. Hence, it
is enough to exclude a neighborhood of the critical latitude
$\beta_{\mbox{\tiny{crit}}}$ which can be taken also exponentially
small. Furthermore the set of points in $\Sset^2$ not covered by
invariant curves of the full system has a measure bounded by the
square root of the perturbation, again exponentially small in
$\eps.$
\end{enumerate}
Summarizing, when arbitrary $g\in SO(3)\setminus{\cal B}$ are
considered the dynamics in $\Sset^2$ is ordered (the points lie on
invariant curves) except for points in a subset of $\Sset^2$ of
exponentially small measure. Furthermore ${\cal B}$ consists of a
neighborhood of the identity of size $O(\eps^{2/3})$ and a
neighborhood of $\beta_{\mbox{\tiny{crit}}}$ which is
exponentially small in $\eps$. \eproof\end{proof}

\bigskip

\begin{cor}\label{C=corol}
For $\eps$ small enough $\Lambda(\eps)<A\eps^3.$
\end{cor}

\begin{proof} It is sufficient to make  remark
\ref{Rem1} more explicit. The differentials of the maps $f_\eps$
increase the length of the vectors in $T\Sset^2$ by a factor of
the form $1+O(\eps)$ and composing with $g\in SO(3)$ produces no
essential changes in the factor. Hence the values of
$\lambda_1(x,g\circ f_\eps)$ are bounded by $C\eps,$ where $C$ is
a positive constant. The contribution to $\Lambda(\eps)$ of the
$g$ to which Theorem \ref{T=final} applies is bounded by $C\eps$
times an exponentially small amount. On the other hand the
contribution of the excluded neighborhood of
$\beta_{\mbox{\tiny{crit}}}$ is also exponentially small.

Therefore the main contribution to $\Lambda(\eps)$ can only come from the
neighborhood of the identity excluded in Theorem \ref{T=final}. But the
Haar
measure of this set is of the order of
$$\int_0^{\eps^{2/3}} (1-\cos(\theta))\dd\theta=O(\eps^2).$$
This bound and the previous one on $\lambda_1(x,g\circ f_\eps)$ give the
result.
\eproof\end{proof}

\bigskip

If we want to consider more ``realistic'' upper bounds it is possible to
proceed along the ideas in remark \ref{Rem2}. A further consideration is
that
the largest stochastic zones are typically associated to the splitting
of separatrices of the hyperbolic fixed points. From \cite{FS} it follows
that
the splitting can be bounded by
$$\exp\left(-\frac{c}{\log(\mu_{\mbox{\tiny{max}}})}\right),$$
where $\mu_{\mbox{\tiny{max}}}$ is the maximal eigenvalue at the fixed
points and $c>0.$ From (\ref{e=maxeigen}) one has that for $\eps$
small $\log(\mu_{\mbox{\tiny{max}}})=\frac{\pi\eps}{2}+O(\eps^2).$ This
``heuristic'' prediction is in good agreement with the observed behavior
for $\eps$ small.

\section{The case of large $\eps$}\label{S=epstoinfty}

Let us call the subbundle of $P\Sset^2$ tangent to the invariant
circles of $f_\eps$ the
{\em horizontal bundle} and denote it by $H$.
As $\eps\to\infty$, a large portion of $P\Sset^2$
is sucked into a small neighborhood of $H$ under
$f_{\eps\#}$.  The measure $\delta_{g(H)}$
on $P\Sset^2$ that is atomic in each tangent space
and supported on $g(H)$ looks more and more like an
invariant measure for $gf_{\eps,\#}$.
These measures integrate to give Lebesgue:
$$m = \int_{g\in SO(3)} \delta_{g(H)} \dd\nu(g).$$
This yields a heuristic argument for why
the inequality in Question~\ref{Q6} should hold
when $\eps = \infty$.  In this section, we make this
argument rigorous by adding some $\delta$-noise
in ${\cal F}_\eps$,
and replacing invariant measures
with stationary measures. We prove:

\begin{theorem}\label{t=converge}
Let $m_{\eps,\delta}$ be defined as in Section 3, and
let $\varphi_{\eps,\delta}$ be the density of $m_{\eps,\delta}$:
$$d m_{\eps,\delta} =  \varphi_{\eps,\delta} dm.$$
There exists $C>0$ such that, for all $\eps,\delta > 0$,
$$\| \varphi_{\eps,\delta} - 1\|_1 < C\delta^{-11}\eps^{-1/2}$$
where $\|\cdot\|_1$ is the $L^1$-norm with respect to Lebesgue
measure $m$ on $P\Sset^2$.
\end{theorem}

This has the corollary:

\begin{cor}\label{c=inftyineq}  There exists a $C>0$ such that for all $\eps,\delta > 0$,
$$|R(\eps,\delta) - R(\eps)| \leq C \delta^{-11} \eps^{-1/2} \log\eps.$$

In particular, for all $\delta>0$,
$$\lim_{\eps\to\infty} |R(\eps,\delta) - R(\eps)| = 0,$$
where $R(\eps,\delta)$ is the random diffused exponent
defined in Section~\ref{S=Heuristics}.
\end{cor}

\bigskip

As we were finishing this paper, unpublished work of
L. Carleson  and T. Spencer came to our attention \cite{CaSp}.
For the standard map on the 2-torus:
$$g_\eps: (x,y)\mapsto (2x + \eps\sin(2\pi x) - y, x)$$
where $\eps$ measures the strength of the
nonlinearity, they prove that by adding a noise of strength
$\exp -\eps^2$ to the element $g_\eps$,
a Lyapunov exponent of order $\log \eps$ can be established.

\bigskip

\noindent{\bf Proof of Corollary~\ref{c=inftyineq}.}
 From the definitions,
\begin{eqnarray*}
R(\eps,\delta) &=& \int_{P\Sset^2} \log\|Tf_\eps v\| \dd m_{\eps,\delta}(v)\\
&=&\int_{P\Sset^2} \log\|Tf_\eps v\| \varphi_{\eps,\delta}(v) \dd m(v),
\end{eqnarray*}
whereas $R(\eps)$ is the integral of $\log\|Tf_\eps v\|$ with respect to $m$.
Hence,
\begin{eqnarray*}
|R(\eps,\delta) - R(\eps)| &=& |\int_{P\Sset^2} \log\|Tf_\eps v\|
(\varphi_{\eps,\delta}(v) -1) \dd m(v)|\\
&\leq& \| \log Tf_\eps\|_\infty  \;\|\;\varphi_{\eps,\delta} -1\;\|_1\\
&\leq & C \delta^{-11}\eps^{-1/2} \log\eps,
\end{eqnarray*}
by Theorem~\ref{t=converge}. \eproof

\bigskip

\noindent{\bf Proof of Theorem~\ref{t=converge}.}

$SO(3)$ acts transitively on $T_1 S^2$ by isometries
and with trivial stabilizer.
 From now on we identify points in $T_1 S^2$ with elements
of $SO(3)$, and use the group structure in writing our formulas.
We will use $x, y, z$ to denote elements of $P\Sset^2$,
$p,q$ for elements of $\Sset^2$, and $(p,v), (q, w)$ for
elements of $P\Sset^2$ (or $T_1\Sset^2$).

Recall from Section 3 that
$$\varphi_{\eps, \delta} = \int_{P\Sset^2} \varphi_{\eps,\delta,g} d\nu(g),$$
where $\varphi_{\eps,\delta,g}$ is the fixed point
of the operator $L_{\eps,\delta,g}$ defined by:
\begin{eqnarray*}
L_{\eps,\delta,g} \psi(x) &=&
\frac{1}{m(B(x,\delta))}\int_{(g {f_\eps}_\sharp)^{-1}B(x,\delta)}\psi(y)\dd
m(y)\\
\end{eqnarray*}

Setting
$$k_{\delta}(x) = \frac{1}{m(B(e,\delta))} 1_{B(e,\delta)}(x),$$
where $e$ is the identity element of $SO(3)$,
we rewrite $L_{\eps,\delta,g}$ as:
\begin{eqnarray}\label{e.idcase}
L_{\eps,\delta,g}\psi (x)
&=& \int_{P\Sset^2} k_{\delta}(x^{-1}g{f_\eps}_\sharp(y)) \psi(y) \dd m(y).
\end{eqnarray}
Let $L_{\eps,\delta} = L_{\eps,\delta,e}$.  It is clear from
(\ref{e.idcase})
that $L_{\eps,\delta,g}\psi(x)  = L_{\eps,\delta} \psi (g^{-1} x)$

Denote by $K_{\eps,\delta}: P\Sset^2\times P\Sset^2\to {\bf R}_+$
the kernel of the operator $L_{\eps,\delta}$, so that
$$K_{\eps,\delta}(x,y) =  k_{\delta}(x^{-1}{f_\eps}_\sharp(y)).$$
Let $\pi: P\Sset^2 \to \Sset^2$ be the projection along tangent fibers.
By averaging along fibers,
we shall approximate $K_{\eps,\delta}$ by a new kernel
$\hat{K}_{\eps,\delta}$ that is constant along fibers
of the second $P\Sset^2$ - factor.
Define $\hat K_{\eps,\delta}: P\Sset^2\times P\Sset^2\to {\bf R}_+$
by
$$\hat K_{\eps,\delta} (x,y) = \int_{\pi^{-1} \pi y}
K_{\eps,\delta}(x, z )\dd m_{\pi y}(z),$$
where, for $p\in \Sset^2$, $m_p$ denotes the disintegration
of $m$ along the fiber $\pi^{-1} p$.
For $g\in SO(3)$, we
obtain a new operator $\hat L_{\eps,\delta,g}$ on
$L^{\infty}(P\Sset^2)$, given by:
$$\hat L_{\eps,\delta,g}\phi (x) = \int_{P\Sset^2}
\hat K_{\eps,\delta}(g^{-1}x,y) \phi(y) \dd m(y).$$
Let $\hat L_{\eps,\delta} = \hat L_{\eps,\delta, e}$.

The next lemma shows that the operators $\hat L_{\eps,\delta,g}$ have
a good averaging property when applied to densities of measures
that project to Lebesgue measure $\mu$ on $\Sset^2$.

\begin{lemma}\label{l=average} Let $\hat K:
P\Sset^2\times P\Sset^2\to {\bf R}_+$ be any
$L^1$ function
such that:
\begin{enumerate}
\item for all $p\in \Sset^2$,
$$\int_{\pi^{-1} p}\int_{P\Sset^2} \hat K(x,y) \dd m_p(x)\dd m(y) = 1;$$
\item if $\pi(y) = \pi(z)$, then for all $x$, $\hat K(x,y) = \hat K(x,z)$.
\end{enumerate}
For $\phi \in L^\infty(P\Sset^2)$, and $g \in SO(3)$,
define $\hat L_g \phi\in  L^\infty(P\Sset^2)$ by
$$\hat L_g \phi(x)  = \int_{P\Sset^2} \hat K(g^{-1} x,y)\phi(y)\dd m(y).$$

Then:
\begin{enumerate}
\item [a.] $\hat \varphi_g = \hat L_g \varphi$
is the unique fixed point of $\hat L_g$, for any
$\varphi\in L^\infty(P\Sset^2)$  that satisfies:
$$\int_{\pi^{-1} p} \varphi(z) \dd m_p(z) = 1,$$
for all $p\in \Sset^2$.
\item[b.] for all $x\in P\Sset^2$, we have:
$$\int_{SO(3)} \hat \varphi_g(x) \dd\nu(g)  = 1.$$

\end{enumerate}
\end{lemma}

\begin{remark} Lemma~\ref{l=average} applies to the operators $\hat
L_{\eps,\delta,g}$. An example of a function $\varphi$ that
satisfies the hypotheses of Lemma~\ref{l=average} is the density
$\varphi_{\eps,\delta,g}$, for any $\eps,\delta,g$.
We do not use that $\hat \varphi_g = \hat L_g \hat \varphi_g $ below.
\end{remark}
\bigskip

\noindent{\bf Proof of Lemma~\ref{l=average}:} Since it is
constant along fibers of the second factor, $\hat K$ projects to a
function on $P\Sset^2\times \Sset^2$, which we shall also call
$\hat K$.

For $g \in SO(3)$, define $\hat \varphi_g$ by:
$$\hat\varphi_g(x)=\hat L_g1(x)=\int_{p\in\Sset^2}\hat K(g^{-1}x,p)\dd\mu(p).$$
We compute directly that, for any $\varphi$ satisfying the
hypotheses of a.,
\begin{eqnarray*}
\hat L_g(\varphi)(x) &=&
\int_{y\in P\Sset^2} \hat K (g^{-1} x, y)\varphi(y) \dd m(y)\\
&=& \int_{p\in \Sset^2}\int_{z\in \pi^{-1} p}
\hat K(g^{-1} x, p)\varphi(z) \dd m_p(z)\dd \mu(p)\\
& = &  \int_{p\in \Sset^2}
K(g^{-1} x, p) \dd \mu(p)\\
& = & \hat \varphi_g(x).
\end{eqnarray*}
To see that $\hat L_g\hat \varphi_g = \hat \varphi_g$ and finish
the proof of a it is now sufficient to verify that
$$\int_{\pi^{-1} p} \hat\varphi_g(z) \dd m_p(z) = 1.$$ But
$\hat\varphi_g(z)= \hat L_g 1 (z)= L_g 1 (z)$ since the function
$1$ is constant and hence constant on fibers. Now $L_g(1)$ is the
density function of a  measure on $P\Sset^2$ which covers Lebesgue
measure on $\Sset^2$ and hence its integral on fibers equals $1$.
This proves a.

Integrating $\hat \varphi_g(x)$ with respect to $g$ we obtain:
\begin{eqnarray*}
\int_{g\in SO(3)} \hat \varphi_g(x)\dd\nu(g)
&= &
\int_{g\in SO(3)} \int_{p\in \Sset^2}
\hat K(g^{-1} x, p) \dd \mu(p) \dd \nu(g)\\
&=&  \int_{y\in P\Sset^2} \int_{p\in \Sset^2}
\hat K(y, p) \dd \mu(p) \dd m(y)\\
& = & 1,
\end{eqnarray*}
completing the proof of b. \eproof

\bigskip

Returning to the proof of Theorem~\ref{t=converge}, let
$\hat \varphi_{\eps, \delta, g} = \hat L_{\eps,\delta,g} 1
= \hat L_{\eps,\delta,g} \varphi_{\eps,\delta,g}$
be the unique fixed point of $\hat L_{\eps,\delta, g}$
given by Lemma~\ref{l=average}.
We now have:
\begin{eqnarray*}
\|\; \varphi_{\eps,\delta} - 1 \; \|_1 & = & \|\; \int_{g\in SO(3)}
(\varphi_{\eps,\delta,g} - 1) \dd\nu(g)\;\|_1\\
&=& \|\; \int_{g\in SO(3)}
(\varphi_{\eps,\delta,g} - \hat\varphi_{\eps,\delta,g}) \dd\nu(g) \;\|_1 +
\|\; \int_{g\in SO(3)} (\hat \varphi_{\eps,\delta,g} - 1) \dd\nu(g)\;\|_1 \\
&\leq & \|\; \int_{g\in SO(3)}
(L_{\eps,\delta, g} \varphi_{\eps,\delta,g }- \hat L_{\eps,\delta, g}
\varphi_{\eps,\delta,g}) \dd\nu(g)\;\|_1 \\
& &+\; \|\;\int_{g\in SO(3)} (\hat \varphi_{\eps,\delta,g } - 1 ) \dd\nu(g)\;\|_1\\
&=& \|\; \int_{g\in SO(3)}
(L_{\eps,\delta} \varphi_{\eps,\delta,g}\circ g^{-1} - \hat
L_{\eps,\delta} \varphi_{\eps,\delta,g}\circ g^{-1}) \dd\nu(g)\;\|_1\\
&\leq & \int_{g\in SO(3)}\|\;
L_{\eps,\delta} \varphi_{\eps,\delta,g}\circ g^{-1} - \hat
L_{\eps,\delta} \varphi_{\eps,\delta,g}\circ g^{-1}\; \|_1 \dd\nu(g)
\end{eqnarray*}
where we used Lemma~\ref{l=average} to obtain the second to last inequality.

Propositions~\ref{p=Ldifference} and \ref{p=linftybound},
which we state and prove below, now imply that
\begin{eqnarray*}
 \int_{g\in SO(3)}\!\!\! \|\; L_{\eps,\delta} \varphi_{\eps,\delta,g}\circ g^{-1} - \hat
L_{\eps,\delta} \varphi_{\eps,\delta,g}\circ g^{-1}\; \|_1  \dd\nu(g)
&\!\!\leq\!\! & C_1\delta^{-3} \eps^{-1/2} \int_{g\in SO(3)} \!\!\!\|\varphi_{\eps,\delta,g}\|_{\infty}
\dd\nu(g)\\
&\!\!\leq\!\! & C_2\delta^{-3} \eps^{-1/2} \int_{g\in SO(3)} \!\!\!\delta^{-8}
\dd\nu(g)\\
&\!\!= \!\!& C \delta^{-11} \eps^{-1/2},
\end{eqnarray*}
completing the proof of the theorem.  It remains to state
and prove Propositions~\ref{p=Ldifference} and \ref{p=linftybound}.

\begin{prop}\label{p=Ldifference} There is a $C>0$ such that,
if $\eps \geq \delta^{-4}$, then for any $\phi\in L^\infty(P\Sset^2)$,
$$\| L_{\eps,\delta}(\phi) - \hat L_{\eps,\delta}(\phi)\|_1
\leq C \|\phi\|_\infty \delta^{-3} \eps^{-1/2}.$$
\end{prop}

\bigskip

\begin{prop}\label{p=linftybound} There exists $C>0$ such that,
for all $\eps>0$ and $\delta>0$,
$$\|\varphi_{\epsilon,\delta,g}\|_\infty \leq C \delta^{-8}.$$
\end{prop}

\bigskip

\noindent{\bf Proof of Proposition~\ref{p=Ldifference}:}
\begin{eqnarray*}
\| L_{\eps,\delta}(\phi) - \hat L_{\eps,\delta}(\phi)\|_1
&\leq &  \|\phi\|_\infty \|K_{\eps,\delta} - \hat K_{\eps,\delta}\|_1 \\
& = &  \|\phi\|_\infty \int |K_{\eps,\delta}(x, y) -
\hat K_{\eps,\delta}(x,y)|\dd m(x)\dd m(y)\\
& = & \|\phi\|_\infty \int |K_{\eps,\delta}(x,y) -
\int_{\pi^{-1} \pi y} K_{\eps, \delta}(x,z) \dd m_{\pi y}(z)
|\dd m(x)\dd m(y).
\end{eqnarray*}

Let $x = (p,v), y=(q,w)$ be elements of $P\Sset^2$.
For a fixed $p,v,q$, the map  $w\mapsto K_{\eps,\delta}((p,v),(q,w))$
is a constant multiple $c\delta^{-3}$ of the characteristic
function for $\pi^{-1} q \; \cap \; {f_{\eps}}_\#^{-1}\; B((p,v),\delta)$.
Note that, for any measurable set $B$ in a probability space with
measure $m$,
the average value of the function $|m(B) - 1_B|$ is
$2 m(B)(1 - m(B))$.
It follows that:
\begin{eqnarray*}
\|\; K_{\eps,\delta}\! -\!
\hat K_{\eps, \delta}\;\|_1 &\!\! =\!\! & 2c\delta^{-3} \int_{(p,v)\in P\Sset^2}
\int_{q\in {f_{\eps}^{-1}} B(p,\delta)}
\! | \; \beta(p,v,q)\;(1\!-\!\beta(p,v,q))
\;|\dd m((p,v)) \dd \mu(q), \\
\end{eqnarray*}
where $$\beta(p,v,q) = m_{q}({f_\eps}_\#^{-1} B((p,v),\delta)).$$

We next show that for some $k>0$, and for $\eps> \delta^{-4}$,
there is a set $G\subset P\Sset^2$, with
$m(G)\geq 1-\eps^{-1/2}$ such that, for all $(p,v)\in G$,
there is a set $G' = G'(p,v) \subset f_{\eps}^{-1} B(p,\delta)$
with $\mu(G') \geq \mu(B(p,\delta)) - \eps^{-1/2}$,
such that, for $q\in G'$,
\begin{eqnarray}\label{e=goodbeta}
\beta(p,v,q) &\leq& k\eps^{-1/2}\quad\hbox{or}\quad \beta(p,v,q) \geq
1-k\eps^{-1/2}.
\end{eqnarray}
This implies that
\begin{eqnarray*}
\|K_{\eps,\delta}\! -\! \hat K_{\eps,\delta}\|_1&\!\!\! \leq\!\!\!
& 2c \delta^{-3}\left( k \eps^{-1/2} \mu(B(p,\delta)) +
2\eps^{-1/2}\right)\\
& \!\!\!\leq\!\!\!  & C\delta^{-3}\eps^{-1/2},
\end{eqnarray*}
which implies the result.

Fix $\delta < 1/2$, and assume that $\eps > \delta^{-4}$.
For $\alpha>0$, denote
by $C_\alpha$ the $\alpha$-neighborhood of the horizontal
bundle $H\subset P\Sset^2$.  In other words, $C_\alpha$ is
the set of lines in $P\Sset^2\setminus \pi^{-1}\{NP,SP\}$ at
angle $\leq \alpha$ with
the latitudinal circles.

It is not difficult
to see that if the distance from $p$ to the poles is
greater than $\eps^{-1/4}$, then
\begin{eqnarray}\label{e=goodzone}
m_{f_\eps^{-1} p} \;(\; f_{\eps\#}^{-1} (C_{\eps^{-1/2}}\cap \pi^{-1}p) \;)\; &\geq&
1- \eps^{-1/2}.
\end{eqnarray}
Let $G_0$ be the set of $p$ for which (\ref{e=goodzone}) holds,
and let $G = \pi^{-1} G_0$.  Then $m(G)\geq 1-\eps^{-1/2}$.

Fix $(p,v)\in G$, and consider a point $f_{\eps} q\in B(p,\delta)\cap G_0$.
The intersection of $B((p,v),\delta)$ with the fiber $\pi ^{-1} f_\eps q$
is an interval $I$. If the endpoints of $I$ are disjoint from
the interval $J = C_{\eps^{-1/2}} \; \cap \; \pi^{-1} f_\eps q$, then
either $I\supset J$ or $I\cap J =\emptyset$. In the former
case, the length of $f_{\eps\#}^{-1} I$ is greater than
the length of $f_{\eps\#}^{-1} J$, which by (\ref{e=goodzone})
is greater than  $1-\eps^{-1/2}$.  In the latter case, the length of
 $f_{\eps\#}^{-1} I$ is less than $\eps^{-1/2}$.
Hence if we let
$G' = G'(p,v)$ be the set of $q\in f_\eps^{-1} B(p,\delta) \cap G_0$
satisfying:
\begin{eqnarray}\label{e=disjoint}
\partial B((p,v),\delta)\cap C_{\eps^{-1/2}} \cap \pi^{-1}f_\eps q
&=&\emptyset,
\end{eqnarray}
then (\ref{e=goodbeta}) holds for all $q\in G'$.
It remains to show that $\mu(G') \geq \mu(B(p,\delta)) -\eps^{-1/2}$.

For $(p,v)\in P\Sset^2$, denote by
$S((p,v),\delta)$ the geodesic sphere of radius $\delta$
centered at $(p,v)$, so $S((p,v),\delta)=\partial B((p,v),\delta)$.
We will use the following
lemma here and later in the proof of Proposition~\ref{p=linftybound}.

\begin{lemma}\label{c=throwaway} There exists $C>0$ such that, for all
$(p,v)\in \pi^{-1} p$, and all $\alpha > 0$,
$$\mu \;\left( \pi \;( S((p,v),\delta )\cap C_\alpha) \right)\leq C\delta\alpha.$$
\end{lemma}

\noindent{\bf Proof.} The claim
follows from the following facts:
\begin{enumerate}
\item On $T_1 \Sset^2\setminus \pi^{-1} \{NP,SP\}$, the subbundle $H$
(regarded as a submanifold) is uniformly transverse to the fibers of
$P \Sset^2$.
\item There exists a $C>0$ such that for all $\delta$ sufficiently
small, and for all $(p,v)\in P\Sset^2$, the intersection
$S((p,v),\delta)\cap H$ is contained
in a smooth curve of length $\leq C \delta$.
\end{enumerate}
The verification of these facts is left as an exercise.\eproof

\bigskip

 From Lemma~\ref{c=throwaway} it follows that:
\begin{eqnarray*}
\mu(G') &\geq& \mu(B(p,\delta)) - \mu(G_0) - \mu(f_\eps^{-1} \pi \;(
S((p,v),\delta )\cap C_{\eps^{-1/2}}) \\
&\geq& \mu(B(p,\delta)) - \eps^{-1/2}
\end{eqnarray*}
This completes the proof of Proposition~\ref{p=Ldifference}.\eproof

\bigskip

\noindent{\bf Proof of Proposition~\ref{p=linftybound}:}
We know that  $\varphi_{\epsilon,\delta,g}$ is a
function in $L^1$ which
satisfies, for all $(p,v) \in P\Sset^2$,
\begin{eqnarray}\label{e=defvarphi}
\varphi_{\epsilon,\delta,g} (p,v) & = & c \delta ^{-3} \int
_{f^{-1}_{\varepsilon \#} B(g^{-1}(p,v),\delta)}
d\varphi_{\epsilon,\delta,g} (q,w) \; \dd \mu(q) \dd m_q(w)
\end{eqnarray}
and, for all $p\in S^2$,
\begin{eqnarray}\label{e=normal}
\int _{\pi^{-1} p} \varphi_{\epsilon,\delta,g} (p,v) \; dm_p(v) & = & 1.
\end{eqnarray}

Then, by (\ref{e=defvarphi}),
the function $\varphi_{\epsilon,\delta,g} $ is continuous,
and therefore has a maximum $M_{\epsilon,\delta,g} $ that we denote by
$M$. The idea is that a H\"older constant for
$\varphi_{\epsilon,\delta,g}$ can be estimated in terms of $M$. Reporting
in (\ref{e=normal})
gives a bound for $M$ which is independent of $\varepsilon , g$.
Since we want to use (\ref{e=normal}) at the end, it suffices to consider the
H\"older constant along the fibers. So, let $(p,v), (p,v') \in
\pi^{-1}p$.  We have:
$$ |\varphi_{\epsilon,\delta,g} (p,v) - \varphi_{\epsilon,\delta,g} (p, v')| \leq c M \delta^{-3} \;
m\; ({f^{-1}_{\varepsilon \#} B(g^{-1}(p,v),\delta)} \;
\Delta
\;{f^{-1}_{\varepsilon \#} B(g^{-1}(p,v'),\delta)} ),$$
where $m$ is Lebesgue measure on $T_1S^2$, and $A\Delta B$
stands for the set of points which belong to only one
of the subsets $A$ or $B$.

\begin{lemma}~\label{l=symdif} There exists $C >0 $ such that, for all
$v,v'\in T_{1,p}\Sset^2$,
$$m\; ({f^{-1}_{\varepsilon \#} B(g^{-1}(p,v),\delta)} \;
\Delta \;{f^{-1}_{\varepsilon \#} B(g^{-1}(p,v'),\delta)} )
\leq C\delta d(v, v')^{1/4}.$$
\end{lemma}

\begin{remark} In Proposition 3.3., we prove that $\varphi_{\epsilon,\delta,g}$
is as smooth as $f_\eps$.  The point of the arguments here is to get
a H{\"o}lder constant independent of $\eps$.
\end{remark}

 From this lemma, it then follows that $M \leq   5 (Cc)^{4} \delta^{-8}/2$, since:
\begin{eqnarray*}
1 & = & \int
_{T_pS^2}
\varphi_{\epsilon,\delta,g} (p,v) \; dv \\
& \geq & M\int _{-\infty }^\infty
(1-Cc\delta^{-2}{|t|^{1/4}})^+ \; dt \\
& = &  \frac{2M}{5} (Cc)^{-4} \delta^{8}.
\end{eqnarray*}

We now prove Lemma~\ref{l=symdif}.

We have the two balls $B((p,v),\delta))$ and
$B((p,v'),\delta)$.  Let $\alpha = \sqrt{d(v,v')}$. We may assume that
$\alpha << \delta$.
The set $ B((p,v),\delta))\Delta B((p,v'),\delta)$ meets
the fiber $\pi^{-1} q$ in a pair of intervals, each of
length $\leq \alpha^2 << \sqrt{\alpha}$.  If the endpoints of
these intervals do not lie in  $C_{\sqrt{\alpha}}$, then
the entire intervals must be disjoint from  $C_{\sqrt{\alpha}}$.
In other words, if
\begin{eqnarray}\label{e=disjoint2}
\pi^{-1} q \cap (S((p,v),\delta)\cup S((p,v'),\delta)) \cap C_{\sqrt{\alpha}} =\emptyset,
\end{eqnarray}
then
\begin{eqnarray}\label{e=disjoint3}
\pi^{-1} q \cap B((p,v),\delta))\Delta B((p,v'),\delta) \cap
C_{\sqrt{\alpha}} = \emptyset.
\end{eqnarray}
Let $G\subset \Sset^2$ be the set of $q$ satisfying (\ref{e=disjoint2}).
By Lemma~\ref{c=throwaway}, $\mu(G)\geq 1 - 2C\delta\sqrt{\alpha}$.

\begin{claim}\label{c=boundregion}
There exists a $C>0$ such that, for all $\alpha \leq 1$,
$p\in \Sset^2$, and $\eps\geq 0$, if $(p,v)\notin C_{\sqrt \alpha}$, then
$$\|T_{(p,v)} f_{\eps,\#}^{-1}\vert_{T_{1,p}\Sset^2}\| \leq C \alpha^{-1}.$$
\end{claim}

\noindent{\bf Proof.} Recall that $\|T_{(p,v)} f_{\eps,\#}^{-1}\vert_{\pi^{-1}p}\| =
\| T_p f_\eps^{-1} v\|^{-2}$.  With respect to the orthonormal basis
of $T_p \Sset^2$ of the form $\{e_1(p), e_2(p)\}$, where
$e_1(p)\in H$ points in the direction of $f_\eps$-twist and
$e_2(p)$ points toward the north pole $NP$, $T_p f_\eps$ takes the form:
$$T_p f_\eps^{-1} = \left(\begin{array}{cc}
1 & -\beta \\
0 & 1
\end{array}
\right),$$
for some $\beta\geq 0$.
A direct computation shows that there exists
a constant $C>0$ such that, for all $\alpha, \beta\geq 0$,
if the angle between a unit vector $v\in {\bf R}^2$ and
the $x$-axis is greater than
$\sqrt{\alpha}$, then: $$\|  \left(\begin{array}{cc}
1 & -\beta \\
0 & 1
\end{array}
\right) v\|^{-2} \leq C \alpha^{-1}$$
 From this the claim follows.\eproof

Claim~\ref{c=boundregion} and (\ref{e=disjoint3}) imply that for $q\!\in\! G$, the derivative of
$f_{\eps\#}$ on \hfill\newline
$\pi^{-1} q\;\cap\;
B((p,v),\delta)\Delta B((p,v'),\delta)$  is bounded:
$$\|T_{(q,w)}f_{\eps\#}\vert_{T_{1,q}\Sset^2}\| \leq \alpha^{-1},$$
for all $w$ such that $(q,w)\in B((p,v),\delta))\Delta B((p,v'),\delta)$.
But for $q\in G$,
\begin{eqnarray*}
m_{f_\eps^{-1}q} \;(\;f_{\eps\#}^{-1}(B((p,v),\delta))\Delta B((p,v'),\delta)\;)
&\leq & \alpha^{-1} m_{q} \;(\;B((p,v),\delta))\Delta B((p,v'),\delta)\;)\\
&\leq& \alpha^{-1} d(v,v')\\
&= & d(v,v')^{1/2}.
\end{eqnarray*}

But then
\begin{eqnarray*}
m\;(\; B((p,v),\delta))\Delta B((p,v'),\delta\;)&\leq&
2C_1 \delta \alpha^{1/2} + C_2 \delta^2 d(v,v')^{1/2}\\
&\leq & C \delta d(v,v')^{1/4},
\end{eqnarray*}
completing the proof of Proposition~\ref{p=linftybound} and of
Theorem~\ref{t=converge}. \eproof

\section{Discussion} \label{S=discussion}

We have wondered \cite {BPSW} about the relationship of the
random Lyapunov exponent of a measure on the space of volume
preserving diffeomorphisms of a manifold to the mean of the
Lyapunov exponents of the individual members.
The point of the question we raised was to
be able to conclude that in a rich enough family of
diffeomorphisms there must be some with positive Lyapunov
exponents, that is to say positive entropy. At question is what
sort of notion of richness would make such a conclusion valid. We
even proposed that much more might conceivably be true, a lower
bound for the mean of the Lyapunov exponents in terms of the
random exponents for orthogonally invariant measures on volume
preserving diffeomorphisms of the sphere. The orthogonal
invariance of the measure was to provide the necessary
``richness".

In the studied family strong numerical evidence has been found
about the existence of such a lower bound when the values of the
stretching parameter $\eps$ are not too small. In some sense
strong stretching has an effect similar to randomization, but it
depends in a clear way on the concrete map. More concretely
\begin{itemize}
\item Even moderate values of $\eps$ like $\eps\ge 10$ are enough to have an
average of the metric entropy larger than the one corresponding to the random
map.
\item There exist unbounded parameters $\eps$ for which islands are born. The
range of existence of these islands is small, but only the islands
associated to fixed points have been considered.
\item For small $\eps$ the estimated average entropy seems positive and
definitely to be much less than the one of the random map. The numerical
evidence is in favor of the existence of exponentially small lower and upper
bounds (in the present example, with an analytic family).
\end{itemize}

The problems in numerically estimating exponents and how to overcome them
have been discussed. A partial analysis of the family of maps has been done for
$\eps$ small. Even a rough estimate of an upper bound of the averaged
entropy is enough to show that the this averaged entropy falls
below any constant multiple of the entropy of the randomized
system, if $\eps$ is sufficiently small.

Finally, the effect of a small randomization of fixed size $\delta$ of the
individual elements of the family ${\mathcal F}_\eps$ is considered. Now the
mean of the local random exponents of the family is indeed asymptotic to the
random exponent of the entire family as $\eps$ tends to infinity; that is,
$R(\eps,\delta)$ and $R(\eps)$ are asymptotic.

\end{document}